\documentclass[9pt]{article}

\usepackage[centertags]{amsmath}
\usepackage{authblk}
\usepackage{multirow}
\usepackage{amssymb}
\usepackage{graphicx}
\usepackage{listing}
\usepackage{mathrsfs}
\usepackage{amsmath,amsthm,amssymb,fontenc,color,textcomp,amsfonts,graphicx,graphics}
\usepackage[colorlinks=true]{hyperref}
\usepackage{hyperref}
\hypersetup{colorlinks = true,linkcolor = black,citecolor = blue}
\usepackage{xcolor}
\usepackage{booktabs}

\theoremstyle{definition}
\newtheorem{defn}{Definition}[section]
\newtheorem{thm}{Theorem}[section]

\newtheorem{Ex}{Example}[section]
\newtheorem{prop}{Proposition}[section]

\usepackage[linesnumbered,ruled,vlined]{algorithm2e}

\usepackage{caption}
\usepackage{subcaption}

\newcommand{\be}{\begin{equation}}
\newcommand{\ee}{\end{equation}}

\begin{document}

\title{Weighted tensor Golub-Kahan-Tikhonov-type methods applied to image processing using
a t-product}


\author{{Lothar Reichel\thanks{\,e-mail: reichel@math.kent.edu} \; and\;
Ugochukwu O. Ugwu\thanks{\,corresponding author e-mail: uugwu@kent.edu}}
\mbox{\ }\\
{\normalsize Department of Mathematical Sciences, Kent State University, OH 44242}}


\maketitle \vspace*{-0.5cm}

\thispagestyle{empty}

\begin{abstract}
\noindent
This paper discusses weighted tensor Golub-Kahan-type bidiagonalization processes
using the t-product. This product was introduced in [M. E. Kilmer and C. D. Martin, 
Factorization strategies for third order tensors, Linear Algebra Appl., 435 (2011), 
pp.~641--658]. A few steps of a bidiagonalization process with a weighted least squares 
norm are carried out to reduce a large-scale linear discrete ill-posed problem to a 
problem of small size. The weights are determined by symmetric positive definite (SPD)
tensors. Tikhonov regularization is applied to the reduced problem. An algorithm for 
tensor Cholesky factorization of SPD tensors is presented. The data is a laterally 
oriented matrix or a general third order tensor. The use of a weighted Frobenius norm 
in the fidelity term of Tikhonov minimization problems is appropriate when the 
noise in the data has a known covariance matrix that is not the identity. We use the discrepancy principle to 
determine both the regularization parameter in Tikhonov regularization and the number of 
bidiagonalization steps. Applications to image and video restoration are considered. 
\vspace{.3cm}

\end{abstract}

\section{Introduction}\label{sec:in}
We are concerned with the iterative solution of large-scale least squares problems of the
form
\begin{equation}
\min_{ \mathcal{\vec{X }}\in\mathbb{R}^{m\times 1\times n}}
\|\mathcal{A} * \mathcal{\vec{X }} - \mathcal{\vec{B}}\|_{\mathcal{M}^{-1}}^2,~~~
\mathcal{A}\in\mathbb{R}^{\ell \times m\times n},~~~
\mathcal{\vec{B}} \in \mathbb{R}^{\ell \times 1\times n},
\label{srhs}
\end{equation}
where $\mathcal{A}=[a_{ijk}]_{i,j,k=1}^{\ell,m,n}\in\mathbb{R}^{\ell\times m\times n}$ is
a third order tensor, $\mathcal{\vec{X}}$ and $\mathcal{\vec{B}}$ are laterally oriented 
$m \times n$ and $\ell \times n$ matrices, respectively, and 
$\|\cdot\|_{\mathcal{M}^{-1}}$ denotes a weighted Frobenius norm; see below. The tensor 
$\mathcal{A}$ specifies the model and the tensor $\mathcal{\vec{B}}$ represents measured 
data that is contaminated by an error $\mathcal{\vec{E}}$, i.e.,
\begin{equation}\label{Btrue}
\mathcal{\vec{B}} = \mathcal{\vec{B}}_\text{true}+\mathcal{\vec{E}},
\end{equation}
where $\mathcal{\vec{B}}_\text{true}\in\mathbb{R}^{\ell\times 1\times n}$ denotes the 
unknown error-free data tensor associated with $\mathcal{\vec{B}}$. The use of a weighted 
Frobenius norm is appropriate when the error $\mathcal{\vec{E}}$ is not white Gaussian. 

We consider minimization problems \eqref{srhs} for which the Frobenius norm of the 
singular tubes of $\mathcal{A}$, which are analogues of the singular values of a matrix, 
decay rapidly to zero, and there are many singular tubes of tiny Frobenius norm of 
different orders of magnitude; see \cite{KM,RU1} for examples of such problems. This 
makes the solution of \eqref{srhs} a linear discrete ill-posed problem. The operator $*$ 
denotes the t-product introduced by Kilmer and Martin \cite{KM}; see Section \ref{sec2} 
for details. 

We will assume that the (unavailable) system of equations
\[
\mathcal{A}*\mathcal{\vec{X}}=\mathcal{\vec{B}}_\text{true}
\]
is consistent and let $\mathcal{\vec{X}}_\text{true}\in\mathbb{R}^{m\times 1\times n}$
denote its unique solution of minimal Frobenius norm. We would like to determine an 
accurate approximation of $\mathcal{\vec{X}}_\text{true}$ given $\mathcal{A}$ and
$\mathcal{\vec{B}}$.

Problems of the form \eqref{srhs} arise, e.g., in image deblurring, see, e.g.,
\cite{KBHH, RU1,RU2,RU3}, where the goal is to recover an accurate approximation of the 
unavailable exact image $\mathcal{\vec{X}}_\text{true}$ by removing blur and noise from an
available degraded image that is represented by $\mathcal{\vec{B}}$. The tensor 
$\mathcal{A}$ represents a blurring operator and typically is ill-conditioned. The 
ill-conditioning of $\mathcal{A}$, and the presence of the error $\mathcal{\vec{E}}$ in 
$\mathcal{\vec{B}}$ make it difficult to solve \eqref{srhs}. In particular, 
straightforward solution of \eqref{srhs} typically gives a meaningless restoration due to 
a large propagated error, which stems from the error $\mathcal{\vec{E}}$ in 
$\mathcal{\vec{B}}$. Regularization, such as by Tikhonov's method, is required to ensure 
that the computed approximate solution of \eqref{srhs} is useful. Regularization 
techniques for \eqref{srhs} when the norm is the standard Frobenius norm have received 
some attention in the literature; see \cite{GIJS,KM,RU1,RU2,RU3}. Here we consider 
replacing the the minimization problem \eqref{srhs} by a penalized weighted least squares 
problem of the form
\begin{equation}\label{4w}
\min_{\mathcal{\vec{X}}\in\mathbb{R}^{m\times 1\times n}}
\left \{\|\mathcal{A*\vec{X}}-\mathcal{\vec{B}}\|^2_{\mathcal{M}^{-1}} +
\mu^{-1}\|\mathcal{\vec{X}}\|^2_{\mathcal{L}^{-1}}\right \},
\end{equation}
where the tensors $\mathcal{M} \in\mathbb{R}^{\ell \times \ell \times n}$ and 
$\mathcal{L} \in\mathbb{R}^{m \times m \times n}$ are symmetric positive definite (SPD) 
with inverses ${\mathcal{M}^{-1}}$ and ${\mathcal{L}^{-1}}$, that determine weighted 
Frobenius norms; see below for a definition. The notion of positive definiteness for 
tensors under the t-product formalism is discussed by Beik et al. \cite{BIJS}. The first 
term in \eqref{4w} is commonly referred to as the \emph{fidelity term} and the second term 
as the \emph{regularization term}.

The Tikhonov regularization problem \eqref{4w} is said to be in \emph{standard form} when 
$\mathcal{M}$ and $\mathcal{L}$ are equal to the identity tensor, which we denote by 
$\mathcal{I}$. It is well known how large-scale problems in standard form can be reduced 
to small problems by using Arnoldi-type or bidiagonalization-type processes; see 
\cite{GIJS,RU1,RU2,RU3}. Tikhonov minimization problems \eqref{4w} with $\mathcal{L}$ or
$\mathcal{M}$ different from $\mathcal{I}$ are said to be in \emph{general form}. The 
choice of the norm in the fidelity term should depend on properties of the error 
$\mathcal{\vec{E}}$ and the choice of the norm in the regularization term should depend on 
properties of the desired solution $\mathcal{X}_\text{true}$. We will discuss these 
choices in Section \ref{sec5}. Other tensor-based solution methods, that do not use the 
t-product, for minimization problems \eqref{4w} in standard form are discussed in 
\cite{BJNR,BNR,GIJB,IGJ}.

Throughout this paper, $\|\mathcal{\vec{X}}\|_{\mathcal{N}}$ denotes the Frobenius norm of
$\mathcal{\vec{X}}$ induced by an SPD tensor 
$\mathcal{N}\in\mathbb{R}^{m\times m\times n}$. We will refer to this norm as 
the $\mathcal{N}$-norm of a tensor column $\mathcal{\vec{X}}$; see Section \ref{sec2} for 
the definition of the $\mathcal{N}$-norm of a general third order tensor 
$\mathcal{X}\in\mathbb{R}^{m\times p\times n}$, $p>1$. The $\mathcal{N}$-norm of 
$\mathcal{\vec{X}}\in\mathbb{R}^{m\times 1\times n}$ is defined as
\[
\|\mathcal{\vec{X}}\|_\mathcal{N} = \sqrt{\big(\mathcal{\vec{X}}^T*\mathcal{N}*
\mathcal{\vec{X}}\big)_{(:,:,1)}}.
\]
Thus, this norm is the square root of the $(1,1,1)$th entry of 
$\mathcal{\vec{X}}^T*\mathcal{N}*\mathcal{\vec{X}}\in\mathbb{R}^{1\times 1\times n}$,
where the superscript $^T$ denotes transposition (defined below).

The quantity $\mu>0$ in \eqref{4w} denotes the regularization parameter. It decides the 
relative influence on the solution of \eqref{4w} of the fidelity and the regularization 
terms. In the present paper, we determine $\mu$ by the discrepancy principle; see, e.g.,
\cite{EHN} for a description and analysis of this principle. Let a bound
\begin{equation}\label{errbd}
\|\mathcal{\vec{E}}\|_{\mathcal{M}^{-1}}\leq\delta
\end{equation}
be known. The discrepancy principle prescribes that $\mu>0$ be determined so that the 
solution $\mathcal{\vec{X}}_\mu$ of \eqref{4w} satisfies
\begin{equation}\label{discr}
\|\mathcal{A}*\mathcal{\vec{X}}_\mu-\mathcal{\vec{B}}\|_{\mathcal{M}^{-1}}=\eta\delta,
\end{equation} 
where $\eta>1$ is a user-specified constant that is independent of $\delta$. Our reason
for using the regularization parameter $\mu$ instead of $1/\mu$ will be commented on in
Section \ref{sec3}. Other techniques, such as generalized cross validation and the L-curve
criterion also can be used to determine the regularization parameter; see, e.g., 
\cite{FRR,Ha,Ki,KR,RR}.

We remark that minimization problems of the form \eqref{4w} also arise in uncertainty 
quantification (UQ) in large-scale Bayesian linear inverse problems; see, e.g., 
\cite{CS,CSBW,SCP}, where $\mathcal{L}$ and $\mathcal{M}$ are suitably defined 
{\it covariance tensors}. The techniques developed in this work can be applied to the 
solution of inverse UQ problems.

It is the purpose of this paper to discuss applications in image and video restoration. We
consider the situation when $\mathcal{M} \neq \mathcal{I}$.  Arridge et al. \cite{ABH} 
observed that for large-scale problems in three space-dimensions, the computation of the 
Cholesky factorization of $\mathcal{L}$ may be prohibitively expensive. We therefore 
derive a solution method that does not require the Cholesky factorization of 
$\mathcal{L}$. 

The normal equations associated with \eqref{4w} are given by
\be \label{normal1}
(\mathcal{A}^T*\mathcal{M}^{-1}*\mathcal{A} + \mu^{-1}\mathcal{L}^{-1})*\mathcal{\vec{X}} 
= \mathcal{A}^T*\mathcal{M}^{-1}*\mathcal{\vec{B}},
\ee
and can be derived analogously as when $\mathcal{L}$ and $\mathcal{M}$ are the identity 
tensor $\mathcal{I}$; see \cite{RU1} for this situation. Since the tensor  $\mathcal{L}$ 
is SPD, equation \eqref{normal1} has a unique solution $\mathcal{X}_\mu$ for any $\mu>0$. 

The change of variables $\mathcal{\vec{X}} = \mathcal{L}*\mathcal{\vec{Y}}$ in
\eqref{normal1} yields the equation
\[
(\mathcal{A}^T*\mathcal{M}^{-1}*\mathcal{A*L} + \mu^{-1}\mathcal{I})*\mathcal{\vec{Y}} =
\mathcal{A}^T*\mathcal{M}^{-1}*\mathcal{\vec{B}}.
\]
Its solution $\mathcal{\vec{Y}}_\mu$ solves the minimization problem
\begin{equation}\label{n4}
\min_{\mathcal{\vec{Y}}\in\mathbb{R}^{m\times 1\times n}}
\{\|\mathcal{A*L*\vec{Y}}-\mathcal{\vec{B}}\|^2_{\mathcal{M}^{-1}} + 
\mu^{-1}\|\mathcal{\vec{Y}}\|^2_\mathcal{L}\}.
\end{equation}
The solution of \eqref{normal1} can be computed as 
$\mathcal{\vec{X}}_\mu = \mathcal{L}*\mathcal{\vec{Y}}_\mu$.

We describe two weighted t-product Golub-Kahan bidiagonalization-type processes for the
approximate solution of \eqref{n4}, and thereby of \eqref{4w}. They generate orthonormal 
tensor bases for the $k$-dimensional tensor Krylov (t-Krylov) subspaces
\begin{equation}\label{krylov}
\mathbb{K}_k(\mathcal{A}^T*\mathcal{M}^{-1}*\mathcal{A*L},\mathcal{A}^T*\mathcal{M}^{-1}*
\mathcal{\vec{B}}) ~~ {\rm and} ~~ \mathbb{K}_k(\mathcal{A*L}*\mathcal{A}^T*\mathcal{M}^{-1},
\mathcal{\vec{B}}),
\end{equation}
where $\mathbb{K}_k(\mathcal{A},\mathcal{\vec{B}}):={\rm span}\{\mathcal{\vec{B}},
\mathcal{A}*\mathcal{\vec{B}},\ldots,\mathcal{A}^{k-1}*\mathcal{\vec{B}}\}$. Each step of
these processes requires two tensor-matrix product evaluations, one with $\mathcal{A}$ and 
one with $\mathcal{A}^T$. Application of a few $k \ll m$ steps of each bidiagonalization 
process to $\mathcal{A}$ reduces the large-scale problem \eqref{n4} to a problem of small 
size. Specifically, $k$ steps of the weighted t-product Golub-Kahan bidiagonalization 
(W-tGKB) process applied to $\mathcal{A}$ reduces $\mathcal{A}$ to a small 
$(k+1) \times k \times n$ lower bidiagonal tensor, while the weighted global tGKB 
(WG-tGKB) process reduces $\mathcal{A}$ to a small $(k+1)\times k$ lower bidiagonal 
matrix. The WG-tGKB process differs from the W-tGKB process in the choice of inner 
product and involves flattening since it reduces equation \eqref{n4} to a problem that
involves a matrix and vectors. Both the W-tGKB and WG-tGKB processes extend the 
generalized Golub-Kahan bidiagonalization process described in \cite{CS} for matrices to
third order tensors. We refer to the solution methods based on the W-tGKB and WG-tGKB 
processes as the weighted t-product Golub-Kahan-Tikhonov (W-tGKT) and the weighted global 
t-product Golub-Kahan-Tikhonov (WG-tGKT) methods, respectively.

We also consider analogues of the minimization problems \eqref{4w} and \eqref{n4} in which
the tensors $\mathcal{\vec{X}}\in\mathbb{R}^{m \times 1\times n}$ and 
$\mathcal{\vec{B}}\in\mathbb{R}^{\ell \times 1\times n}$ are replaced by general third 
order tensors $\mathcal{X}\in\mathbb{R}^{m \times p\times n}$ and 
$\mathcal{B} \in\mathbb{R}^{\ell \times p\times n}$, respectively, for some $p>1$. That 
is, we consider minimization problems of the form
\begin{equation}
\min_{ \mathcal{X}\in\mathbb{R}^{m\times p\times n}}
\|\mathcal{A} * \mathcal{X} - \mathcal{B}\|_{\mathcal{M}^{-1}}^2.
\label{3tensor}
\end{equation}
Problems of this kind arise when restoring a color image or a sequence of consecutive 
gray-scale video frames. For instance, a color
image restoration problem with RGB channels corresponds to $p=3$, whereas when restoring a
gray-scale video, $p$ is the number of consecutive video frames. The minimization problems 
analogous to \eqref{4w} and \eqref{n4} are given by
\begin{equation}\label{4m}
\min_{\mathcal{X}\in\mathbb{R}^{m\times p\times n}}
\{\|\mathcal{A*X}-\mathcal{B}\|^2_{\mathcal{M}^{-1}} + 
\mu^{-1}\|\mathcal{X}\|^2_{\mathcal{L}^{-1}}\}
\end{equation}
and 
\begin{equation}\label{nn4}
\min_{\mathcal{Y}\in\mathbb{R}^{m\times p\times n}}
\{\|\mathcal{A*L*Y}-\mathcal{B}\|^2_{\mathcal{M}^{-1}} + 
\mu^{-1}\|\mathcal{Y}\|^2_{\mathcal{L}}\},
\end{equation}
respectively.

We present three solution methods for \eqref{4m}. The first two of them work with each 
tensor column $\mathcal{\vec{B}}_j$, $j=1,2,\dots,p$, of $\mathcal{B}$ independently, 
i.e., they consider \eqref{nn4} as $p$ separate minimization problems and apply the W-tGKT 
or WG-tGKT methods to each one of these problems. The resulting methods for \eqref{4m} are
referred to as the W-tGKT$_p$ and WG-tGKT$_p$ methods, respectively, and sometimes simply 
as $p$-methods. The third method works with the lateral slices of $\mathcal{B}$ simultaneously
and will be referred to as the weighted generalized global tGKT (WGG-tGKT) method. This 
method uses the weighted generalized global tGKB (WGG-tGKB) process introduced in Section
\ref{sec4} to reduce \eqref{nn4} to a problem of small size. The WGG-tGKT method requires 
less CPU time than the $p$-methods because it uses larger chunks of data at a time than 
the other methods in our comparison.

Numerical experiments in \cite{RU1, RU2, RU3} and in Section \ref{sec5} show the 
$p$-methods to often yield restorations of higher quality than the WGG-tGKT method. This 
behavior may be attributed to the facts that the $p$-methods determine $p$ regularization
parameters, and generally use $p$ different t-Krylov subspaces \eqref{krylov}. The 
WGG-tGKT and WG-tGKT$_p$ methods involve flattening and require additional product 
definition to the t-product. They are related to the global tensor Krylov subspace methods
described in \cite{BJNR,BNR,GIJB,GIJS,IGJ,RU1,RU2,RU3}.

The organization of this paper is as follows. Section \ref{sec2} introduces notation and 
preliminaries associated with the t-product, while Section \ref{sec3} discusses the W-tGKT
and W-tGKT$_p$ methods as well as the W-tGKB process. Section \ref{sec4} presents the 
weighted global t-Krylov subspace methods. The WGG-tGKT method and WGG-tGKB process are 
described in Section \ref{sec4.1}. Section \ref{sec4.2} presents the WG-tGKB process and
the WG-tGKT and WG-tGKT$_p$ methods. In the computed examples in Section \ref{sec5}, we 
take $\mathcal{M \neq I}$ and discuss the performance of these methods when 
$\mathcal{L \neq I}$ and $\mathcal{L = I}$. Concluding remarks are presented in Section 
\ref{sec6}.

\section{Notation and Preliminaries}\label{sec2}
This section reviews and modifies results by El Guide et al. \cite{GIJS} and Kilmer et al.
\cite{KBHH} to suit the current framework. In this paper, a tensor is a multidimensional 
array of real numbers. We denote third order tensors by calligraphic script, e.g., 
$\mathcal{A}$, matrices by capital letters, e.g., $A$, vectors by lower case letters, 
e.g., $a$, and tubal scalars (tubes) by boldface letters, e.g., $\mathbf{a}$. Tubes of 
third order tensors are obtained by fixing any two of the indices, while slices are 
obtained by keeping one of the indices fixed; see Kolda and Bader \cite{KB}. The $j$th 
lateral slice (tensor column) is a laterally oriented matrix denoted by 
$\mathcal{\vec{A}}_j$. The $k$th frontal slice is denoted by $\mathcal{A}^{(k)}$ and is a 
matrix.

Given a third order tensor $\mathcal{A} \in \mathbb{R}^{\ell \times m \times n}$ with
$\ell \times m$ frontal slices $\mathcal{A}^{(i)}$, $i=1,2,\ldots,n$, the operator
$\mathtt{unfold}(\mathcal{A})$ returns an $\ell n\times m$ matrix with the frontal slices,
whereas the $\mathtt{fold}$ operator folds back the unfolded tensor $\mathcal{A}$, i.e.,
\begin{equation*}
\mathtt{unfold}(\mathcal{A}) = \begin{bmatrix}
\mathcal{A}^{(1)}\\
\mathcal{A}^{(2)}\\
\vdots\\
\mathcal{A}^{(n)}
\end{bmatrix}, \;\;\;\;\; \mathtt{fold(unfold(\mathcal{A})) = \mathcal{A}}.
\end{equation*}
The operator $\mathtt{bcirc}(\mathcal{A})$ generates an $\ell n \times mn$ block-circulant
matrix with $\mathtt{unfold}(\mathcal{A})$ forming the first block column,
\begin{equation*}
\mathtt{bcirc}(\mathcal{A}) = \begin{bmatrix}
\mathcal{A}^{(1)} & \mathcal{A}^{(n)} & \dots & \mathcal{A}^{(2)}\\
\mathcal{A}^{(2)} & \mathcal{A}^{(1)} & \dots & \mathcal{A}^{(3)}\\
\vdots & \vdots & \ddots & \vdots\\
\mathcal{A}^{(n)} & \mathcal{A}^{(n-1)} & \dots & \mathcal{A}^{(1)}
\end{bmatrix}.
\end{equation*}

\begin{defn}{(t-product \cite{KM})}\label{def: 2.1}
Let $\mathcal{A}\in\mathbb{R}^{\ell\times m\times n}$ and
$\mathcal{B}\in\mathbb{R}^{m\times p\times n}$. Then the t-product $\mathcal{A*B}$ is the
tensor $\mathcal{C}\in\mathbb{R}^{\ell\times p\times n}$ defined by
\begin{equation}\label{tenA}
\mathcal{C}=\mathtt{fold}(\mathtt{bcirc}(\mathcal{A}) \cdot \mathtt{unfold}(\mathcal{B})),
\end{equation}
where $\cdot$ stands for the standard matrix-matrix product.
\end{defn}

The block circulant matrix $\mathtt{bcirc}(\mathcal{A})$ can be block diagonalized by the 
discrete Fourier transform (DFT) matrix combined with a Kronecker product. Suppose 
$\mathcal{A}\in\mathbb{R}^{\ell\times m\times n}$ and let $F_n$ be an $n\times n$ unitary 
DFT matrix defined by
\begin{equation*}
F_n = \frac{1}{\sqrt{n}}\begin{bmatrix}
1 & 1 & 1 & \dots & 1\\
1 & \omega & \omega^2 & \dots & \omega^{n-1}\\
1 & \omega^2 & \omega^4 & \dots & \omega^{2(n-1)}\\
\vdots & \vdots & \vdots & \ddots & \vdots\\
1 & \omega^{n-1} & \omega^{n-1} & \dots & \omega^{(n-1)(n-1)}\\
\end{bmatrix}
\end{equation*}
with $\omega = e^{\frac{-2\pi i}{n}}$ and $i^2 = -1$. Then
\begin{equation}\label{flop}
\bar{A}:=\mathtt{blockdiag}(\widehat{A}^{(1)},\widehat{A}^{(2)},\dots,
\widehat{A}^{(n)})= (F_n \otimes I_\ell)\cdot \mathtt{bcirc}(\mathcal{A})\cdot
(F_n^* \otimes I_m),
\end{equation}
where $\otimes$ stands for Kronecker product and $F_n^*$ denotes the conjugate transpose 
of $F_n$. The matrix $\bar{A}\in\mathbb{R}^{\ell n\times mn}$ is block diagonal with the 
diagonal blocks $\widehat{A}^{(i)}\in\mathbb{R}^{\ell\times m}$, $i=1,2,\dots,n$. The 
matrices $\widehat{A}^{(i)}$ are the frontal slices of the tensor $\widehat{\mathcal{A}}$
obtained by applying the FFT along each tube of $\mathcal{A}$. Each matrix 
$\widehat{A}^{(i)}$ may be dense and have complex entries unless certain symmetry 
conditions hold; see \cite{KM} for further details. Throughout this paper, we often will 
denote objects that are obtained by taking the FFT along the third dimension with a 
widehat over the argument, i.e., $\widehat{\cdot}$~.

Using \eqref{flop}, the t-product \eqref{tenA} can be evaluated as
\begin{equation}\label{tprodfft}
\mathcal{A*B} = \mathtt{fold}\Big((F_n^* \otimes I_\ell)\cdot
\big((F_n \otimes I_\ell)\cdot\mathtt{bcirc}(\mathcal{A})\cdot(F_n^* \otimes I_m)\big)\cdot
(F_n \otimes I_m)\cdot \mathtt{unfold}(\mathcal{B})\Big).
\end{equation}
It is easily shown that by taking the FFT along the tubes of
$\mathcal{B}\in\mathbb{R}^{m\times p\times n}$, we can compute
$(F_n \otimes I_m)\cdot \mathtt{unfold}(\mathcal{B})$ in $\mathcal{O}(pmn\log_2(n))$ 
arithmetic floating point operations (flops); see Kilmer and Martin \cite{KM} for details.

The computations \eqref{tprodfft} can be easily implemented in MATLAB. Using MATLAB
notation, let $\mathcal{\widehat{C}}:=\mathtt{fft}(\mathcal{C},[\;],3)$ be the tensor
obtained by applying the FFT along the third dimension and let 
$\mathcal{\widehat{C}}(:,:,i)$ denote the $i$th frontal slice of $\mathcal{\widehat{C}}$.
Then the t-product \eqref{tprodfft} can be computed by taking the FFT along the tubes of
$\mathcal{A}$ and $\mathcal{B}$ to get 
$\mathcal{\widehat{A}} = \mathtt{fft}(\mathcal{A},[\;],3)$ and
$\mathcal{\widehat{B}} = \mathtt{fft}(\mathcal{B},[\;],3)$, followed by a matrix-matrix 
product of each
pair of the frontal slices of $\mathcal{\widehat{A}}$ and $\mathcal{\widehat{B}}$, i.e.,
\begin{equation*}
\mathcal{\widehat{C}}(:,:,i) = \mathcal{\widehat{A}}(:,:,i)\cdot
\mathcal{\widehat{B}}(:,:,i),\quad i = 1, 2, \dots, n,
\end{equation*}
and then taking the inverse FFT along the third
dimension to obtain $\mathcal{C} = \mathtt{ifft}(\mathcal{\widehat{C}},[\;],3)$.

Let $\mathcal{A} \in \mathbb{R}^{\ell \times m \times n}$. Then the tensor transpose
$\mathcal{A}^T\in\mathbb{R}^{m \times \ell \times n}$ is the tensor obtained by first
transposing each one of the frontal slices of $\mathcal{A}$, and then reversing the order
of the transposed frontal slices 2 through $n$; see \cite{KM}. The tensor transpose is 
analogous to the matrix transpose. Specifically, if $\mathcal{A}$ and
$\mathcal{B}$ are two tensors such that $\mathcal{A*B}$ and $\mathcal{B}^T*\mathcal{A}^T$
are defined, then $(\mathcal{A*B})^T = \mathcal{B}^T*\mathcal{A}^T$. The identity 
$\mathcal{I} \in \mathbb{R}^{m \times m \times n}$ is a tensor whose first frontal slice, 
$\mathcal{I}^{(1)}$, is the $m \times m$  identity matrix, and all other frontal slices, 
$\mathcal{I}^{(i)}$, $i = 2,3, \dots, n$, are zero matrices; see \cite{KM}.

The notion of orthogonality under the t-product formalism is well defined and similar to 
the matrix case; see Kilmer and Martin \cite{KM}. A tensor 
$\mathcal{Q}\in\mathbb{R}^{m\times m\times n}$ is said to be orthogonal if 
$\mathcal{Q}^T*\mathcal{Q} = \mathcal{Q}*\mathcal{Q}^T = \mathcal{I}$. The lateral slices 
of $\mathcal{Q}$ are orthonormal if they satisfy
\begin{equation*}
\mathcal{Q}^T(:,i,:)*\mathcal{Q}(:,j,:) = \left\{
\begin{array}{ll}
{\bf e}_1 & i=j,\\
{\bf 0} & i \neq j,
\end{array}
\right.
\end{equation*}
where $\mathbf{e}_1\in\mathbb{R}^{1\times 1\times n}$ is a tubal scalar with the
$(1,1,1)$th entry equal to $1$ and the remaining entries zero. Given the tensor 
$\mathcal{A}$ and an orthogonal tensor $\mathcal{Q}$ of compatible dimension, we have
\[
\| \mathcal{Q*A}\|_F = \|\mathcal{A}\|_F;
\]
see \cite{KM} for a proof. The notion of \textit{partial orthogonality} is similar as in 
matrix theory. The tensor $\mathcal{Q}\in\mathbb{R}^{\ell\times m\times n}$ with $\ell>m$
is said to be partially orthogonal if $\mathcal{Q}^T*\mathcal{Q}$ is well defined and 
equal to the identity tensor $\mathcal{I}\in\mathbb{R}^{m\times m\times n}$; see 
\cite{KM}.

A tensor $\mathcal{A}\in\mathbb{R}^{m\times m\times n}$ has an inverse, denoted by 
$\mathcal{A}^{-1}$, provided that $\mathcal{A}*\mathcal{A}^{-1}=\mathcal{I}$ and 
$\mathcal{A}^{-1}*\mathcal{A}=\mathcal{I}$, whereas a tensor is said to be f-diagonal if 
each frontal slice of the tensor is a diagonal matrix; see \cite{KM}.

Let $\mathcal{A}\in\mathbb{R}^{\ell\times m\times n}$, $\ell \geq m$. Then the tensor 
singular value decomposition (tSVD), introduced by Kilmer and Martin \cite{KM}, is given 
by
\[
\mathcal{A} = \mathcal{U}*\mathcal{S}*\mathcal{V}^T,
\]
where $\mathcal{U}\in\mathbb{R}^{\ell\times\ell\times n}$ and
$\mathcal{V}\in\mathbb{R}^{m\times m\times n}$ are orthogonal tensors, and the tensor
\[
\mathcal{S}={\rm diag}[\mathbf{s}_1,\mathbf{s}_2,\dots,\mathbf{s}_m]\in
\mathbb{R}^{\ell\times m\times n}
\]
is f-diagonal with singular tubes ${\mathbf{s}}_j\in\mathbb{R}^{1\times 1\times n}$,
$j =1,2,\dots,m$, ordered according to
\[
\|\mathbf{s}_1\|_F\geq\|\mathbf{s}_2\|_F\geq\cdots\geq\|\mathbf{s}_m\|_F.
\]
The tubal rank of $\mathcal{A}$ is the number of nonzero singular tubes of $\mathcal{A}$;
see Kilmer et al. \cite{KBHH}.

The Frobenius norm of a tensor column $\mathcal{\vec{X}}\in\mathbb{R}^{m\times 1\times n}$
is given by
\begin{equation}\label{defn: 2.7}
\|\mathcal{\vec{X}}\|_F = \sqrt{\big(\mathcal{\vec{X}}^T*\mathcal{\vec{X}}\big)_{(:,:,1)}};
\end{equation}
see \cite{KBHH}. Thus, the square of the Frobenius norm of $\mathcal{\vec{X}}$ is the
$(1,1,1)$th entry of the tubal scalar 
$\mathcal{\vec{X}}^T*\mathcal{\vec{X}}\in\mathbb{R}^{1\times 1\times n}$, which is denoted
by $(\mathcal{\vec{X}}^T*\mathcal{\vec{X}}\big)_{(:,:,1)}$.

A tensor $\mathcal{M}\in\mathbb{R}^{m\times m\times n}$ is said to be positive definite 
(PD) if it satisfies 
\[
\big(\mathcal{\vec{X}}^T*\mathcal{M}*\mathcal{\vec{X}}\big)_{(:,:,1)} > 0
\]
for all nonzero tensor columns $\mathcal{\vec{X}}\in\mathbb{R}^{m\times 1\times n}$;
see \cite{BIJS}. Moreover, $\mathcal{M}$ is symmetric if each frontal slice
$\mathcal{\widehat{M}}^{(i)}$, $i = 1,2,\dots, n$, of $\mathcal{\widehat{M}}$ is 
Hermitian. Hence, a tensor $\mathcal{M}$ is SPD if $\mathcal{\widehat{M}}^{(i)}$ is 
Hermitian PD; see \cite{KBHH}.

Let $\mathcal{M}\in\mathbb{R}^{m\times m\times n}$ be an SPD tensor. The 
$\mathcal{M}$-norm of a general third order tensor 
$\mathcal{X}\in\mathbb{R}^{m\times p\times n}$ with $p>1$ is given by
\[
\|\mathcal{X}\|_\mathcal{M} = 
\sqrt{{\tt trace}\Big(\big(\mathcal{X}^T*\mathcal{M}*\mathcal{X}\big)_{(:,:,1)}\Big)}.
\]
The quantity $\big(\mathcal{X}^T*\mathcal{M}*\mathcal{X}\big)_{(:,:,1)}$ is a $p \times p$
matrix. It represents the first frontal slice of the tensor 
$\mathcal{X}^T*\mathcal{M}*\mathcal{X}\in\mathbb{R}^{p\times p\times n}$, and 
${\tt trace}(M)$ denotes the trace of the matrix $M$.

Algorithm \ref{Alg: 01} takes a nonzero tensor
$\mathcal{\vec{X}} \in \mathbb{R}^{m\times 1\times n}$ and returns a normalized tensor
$\mathcal{\vec{V}}\in\mathbb{R}^{m\times 1\times n}$ and a tubal scalar
$\mathbf{a}\in\mathbb{R}^{1\times 1\times n}$, such that
\begin{equation*}
\mathcal{\vec{X}} = \mathcal{\vec{V}}* \mathbf{a} ~~{\rm and}~~
\| \mathcal{\vec{V}}\|_{\mathcal{M}} = 1.
\end{equation*}
The tubal scalar $\mathbf{a}$ is not necessarily invertible; see \cite{KBHH}. We remark
that $\mathbf{a}$ is invertible only if there is a tubal scalar $\mathbf{b}$ such that 
$\mathbf{a*b} = \mathbf{b*a} = \mathbf{e}_1$. The scalar $\mathbf{a}^{(j)}$ is the $j$th 
face of the $1\times 1\times n$ tubal scalar $\mathbf{a}$, while $\mathcal{\vec{V}}^{(j)}$
is a vector with $m$ entries, and the $j$th face of 
$\mathcal{\vec{V}}\in\mathbb{R}^{m\times 1\times n}$. Steps 3 and 7 of Algorithm 
\ref{Alg: 01} use the fact that for an SPD matrix $M \in \mathbb{R}^{m \times m}$ and 
$x \in \mathbb{R}^{m}$, we have
\[
\| x\|_M = \sqrt{x^TMx}.
\]
An analogue of Algorithm \ref{Alg: 01} for $\mathcal{M=I}$ is described in \cite{KBHH}.

\vspace{.3cm}
\begin{algorithm}[H]
\SetAlgoLined
\KwIn{$\mathcal{\vec{X}} \in \mathbb{R}^{m \times 1 \times n} \neq \mathcal{\vec{O}}$, $\mathcal{M}$ is an $m \times m\times n$ SPD tensor}
\KwOut{ $\mathcal{\vec{V}}$, $\mathbf{a}$ with $\mathcal{\vec{X}}=\mathcal{\vec{V}}*\mathbf{a}$ and $\| \mathcal{\vec{V}}\|_{\mathcal{M}} = 1$}
$\mathcal{\vec{V}} \leftarrow \mathtt{fft}(\mathcal{\vec{X}},[\;],3)$, $\mathcal{M} \leftarrow \mathtt{fft}(\mathcal{M},[\;],3)$\\
\For{$j = 1,2,\dots,n$}{
$\mathbf{a}^{(j)} \gets \|\mathcal{\vec{V}}^{(j)}\|_{\mathcal{M}^{(i)}}\;\;\;$ ($\mathcal{\vec{V}}^{(j)}$ is a vector)\\
\eIf{$\mathbf{a}^{(j)} > \mathtt{tol}$}{
$\mathcal{\vec{V}}^{(j)} \gets \frac{1}{\mathbf{a}^{(j)}} \mathcal{\vec{V}}^{(j)}$\\
}{
$\mathcal{\vec{V}}^{(j)} \gets \mathtt{randn}(m,1); \;\; \mathbf{a}^{(j)} \gets \|\mathcal{\vec{V}}^{(j)}\|_{\mathcal{M}^{(i)}};\;\;
\mathcal{\vec{V}}^{(j)} \gets \frac{1}{\mathbf{a}^{(j)}} \mathcal{\vec{V}}^{(j)}; \;\; \mathbf{a}^{(j)} \gets 0$\\
} }
$\mathcal{\vec{V}} \gets \mathtt{ifft}(\mathcal{\vec{V}},[\;],3); \;\; \mathbf{a} \gets \mathtt{ifft}(\mathbf{a},[\;],3)$
\caption{Normalize($\mathcal{\vec{X}},\mathcal{M}$)}
\label{Alg: 01}
\end{algorithm}\vspace{.3cm}

Algorithm \ref{cholesky} determines the tensor Cholesky decomposition of an SPD tensor 
$\mathcal{M}$. This algorithm is mimetic of the Cholesky factorization of a matrix $M$.
Step $3$ of Algorithm \ref{cholesky} carries out the Cholesky factorization of the 
Hermitian matrix $\mathcal{\widehat{M}}^{(i)}$ in the Fourier domain. This algorithm will
be used in Section \ref{sec5}.

\vspace{.3cm}
\begin{algorithm}[H]
\SetAlgoLined
\KwIn{SPD tensor $\mathcal{M}\in\mathbb{R}^{m \times m\times n}$}
\KwOut{Tensor Cholesky factor $\mathcal{R}\in\mathbb{R}^{m\times m\times n}$, where 
$\mathcal{M} = \mathcal{R}^T*\mathcal{R}$}
$\mathcal{\widehat{M}} = {\tt fft}(\mathcal{M},[\;],3)$\\
\For{$i = 1$ \bf{to} $n$}{
$\mathcal{\widehat{R}}^{(i)} \leftarrow {\tt chol}(\mathcal{\widehat{M}}^{(i)})$, 
where {\tt chol} denotes MATLAB's Cholesky factorization operator. Thus, 
$\mathcal{\widehat{R}}^{(i)}$ is an upper triangular complex matrix.
}
$\mathcal{R} \gets {\tt ifft}(\widehat{\mathcal{R}},[\;],3)$
\caption{Tensor Cholesky factorization}
\label{cholesky}
\end{algorithm}\vspace{.3cm}

Introduce the tensors
\[
\mathbb{C}_k:=[\mathcal{C}_1,\mathcal{C}_2,\dots,\mathcal{C}_k]\in
\mathbb{R}^{m\times kp\times n}, \;\;\;\;\;
\mathcal{C}_k:=[\mathcal{\vec{C}}_1, \mathcal{\vec{C}}_2,\dots,\mathcal{\vec{C}}_k]\in
\mathbb{R}^{m\times k\times n},
\]
where $\mathcal{C}_j\in\mathbb{R}^{m\times p\times n}$ and
$\mathcal{\vec{C}}_j\in\mathbb{R}^{m\times 1\times n}$, $j = 1,2,\dots,k$. Suppose that
$y = [y_1,y_2, \dots, y_k]^T\in\mathbb{R}^k$. El Guide et al. \cite{GIJS} define two 
products denoted by $\circledast$,
\[
\mathbb{C}_k\circledast y=\sum_{j=1}^k y_j\mathcal{C}_j, \;\; \;\;\;
\mathcal{C}_k\circledast y=\sum_{j=1}^k y_j \mathcal{\vec{C}}_j.
\]

We conclude this section by modifying some notions introduced by El Guide et al. 
\cite{GIJS}. Let $\mathcal{M}=[m_{iik}]\in \mathbb{R}^{m\times m\times n}$ and consider 
the tensors 
$\mathcal{C}=[c_{ijk}],\mathcal{D}=[d_{ijk}]\in \mathbb{R}^{m\times p\times n}$ with 
lateral slices $\mathcal{\vec{C}}=[c_{i1k}], \mathcal{\vec{D}}=[d_{i1k}]\in
\mathbb{R}^{m\times 1\times n}$. Introduce the inner products
\begin{equation*}
\langle\mathcal{C},\mathcal{D}\rangle=
\sum_{i=1}^m\sum_{j=1}^p\sum_{k=1}^nc_{ijk}d_{ijk},
\;\;\;\;\;\; \langle\mathcal{\vec{C}},\mathcal{\vec{D}}\rangle=
\sum_{i=1}^m\sum_{j=1}^p\sum_{k=1}^n c_{i1k}d_{i1k}.
\end{equation*}
Then
\begin{equation*}
\langle\mathcal{C},\mathcal{D}\rangle_\mathcal{M}=
\sum_{i=1}^m\sum_{j=1}^p\sum_{k=1}^nc_{ijk}[\mathcal{M*D}]_{ijk},
\;\;\;\;\;\; \langle\mathcal{\vec{C}},\mathcal{\vec{D}}\rangle_\mathcal{M}=
\sum_{i=1}^m\sum_{j=1}^p\sum_{k=1}^n c_{i1k}[\mathcal{M*\vec{D}}]_{i1k}.
\end{equation*}
Specifically, $\langle\mathcal{C},\mathcal{D}\rangle_\mathcal{M} = \langle\mathcal{C},\mathcal{M*D}\rangle$, 
and $\langle\mathcal{C},\mathcal{\vec{D}}\rangle_\mathcal{M} = \langle\mathcal{C},\mathcal{M*\vec{D}}\rangle$.
Let
\begin{equation*}
\begin{split}
\mathbb{A}:=[\mathcal{A}_1,\mathcal{A}_2,\dots,\mathcal{A}_m]\in
\mathbb{R}^{\ell \times sm \times n}, \;\;\;\;\;
\mathbb{B}:=[\mathcal{B}_1,\mathcal{B}_2,\dots,\mathcal{B}_p]\in
\mathbb{R}^{\ell \times sp \times n},\\
\mathcal{A}:=[\mathcal{\vec{A}}_1,\mathcal{\vec{A}}_2,\dots,\mathcal{\vec{A}}_m]\in
\mathbb{R}^{\ell \times m \times n}, \;\;\;\;\;
\mathcal{B}:=[\mathcal{\vec{B}}_1,\mathcal{\vec{B}}_2,\dots,\mathcal{\vec{B}}_p]\in
\mathbb{R}^{\ell \times p \times n},
\end{split}
\end{equation*}
where $\mathcal{A}_i\in\mathbb{R}^{\ell\times s\times n}$,
$\mathcal{\vec{A}}_i\in\mathbb{R}^{\ell\times 1\times n}$, $i = 1,2,\dots,m$, and
$\mathcal{B}_j\in\mathbb{R}^{\ell\times s\times n}$,
$\mathcal{\vec{B}}_j\in\mathbb{R}^{\ell\times 1\times n}$ $j = 1,2,\dots,p$. The weighted 
T-diamond products $\mathbb{A}^T \Diamond_\mathcal{M} \mathbb{B}$ and 
$\mathcal{A}^T \Diamond_\mathcal{M} \mathcal{B}$ yield $m\times p$ matrices given by
\begin{equation*}
[\mathbb{A}^T\Diamond_\mathcal{M}\mathbb{B}]_{ij}=\langle \mathcal{A}_i, \; 
\mathcal{B}_j\rangle_\mathcal{M}, \;\;\;\;\;
[\mathcal{A}^T\Diamond_\mathcal{M}\mathcal{B}]_{ij}=\langle \mathcal{\vec{A}}_i, \;
\mathcal{\vec{B}}_j\rangle_\mathcal{M}, \;\;\; i = 1,2,\dots,m, \;\; j= 1,2,\dots,p.
\end{equation*}

\section{The Weighted t-Product Golub-Kahan-Tikhonov Method}\label{sec3}
This section extends the generalized Golub-Kahan bidiagonalization ({\tt gen}$\_$GKB) 
process described in \cite{CS} for matrices to third order tensors. The {\tt gen}$\_$GKB 
process was first proposed by Benbow \cite{B} for generalized least squares problems. 
Applications of this process to uncertainty quantification in large-scale Bayesian linear
inverse problems have recently been described in \cite{CS,CSBW,SCP}. For additional 
applications, see, e.g., \cite{A, AO, OA}. We will refer to the tensor version of the 
{\tt gen}$\_$GKB process as the weighted t-product Golub-Kahan bidiagonalization (W-tGKB) 
process. The latter process is described by Algorithm \ref{Alg: 3} below. Its global 
versions are presented in Section \ref{sec4}.

We also discuss the computation of an approximate solution of the minimization problems 
\eqref{4w} and \eqref{4m} with the aid of the partial W-tGKB process. Typically, only a 
few, say $k \ll m$, steps of this process are required to reduce large-scale problems 
\eqref{4w} and \eqref{4m} to problems of small size. Specifically, $k$ steps of the W-tGKB
process applied to $\mathcal{A}$ in \eqref{4w} with initial tensor $\mathcal{\vec{B}}$ 
reduces $\mathcal{A}$ to a small lower bidiagonal tensor, which we  denote by 
$\mathcal{\bar{P}}_k \in \mathbb{R}^{{(k+1)}\times k\times n}$. The reduction process is 
described by Algorithm \ref{Alg: 3}.

\vspace{.3cm}
\begin{algorithm}[H]
\SetAlgoLined
\KwIn{$ \mathcal{A} \in \mathbb{R}^{\ell \times m \times n}$, $\mathcal{B} \in
\mathbb{R}^{\ell \times p \times n}$ such that 
$\mathcal{A}^T*\mathcal{M}^{-1}*\mathcal{B}\neq\mathcal{O}$, 
$\mathcal{L} \in \mathbb{R}^{m \times m \times n}$, $\mathcal{M} \in 
\mathbb{R}^{\ell \times \ell \times n}$ are SPD, $k\geq 1$}
$[\mathcal{\vec{Q}}_1, {\mathbf{z}}_1] \leftarrow \mathtt{Normalize}(\mathcal{\vec{B}}, 
\mathcal{M}^{-1})$ with $\mathbf{z}_1$ invertible\\
$[\mathcal{\vec{W}}_1, {\mathbf{c}}_1] \leftarrow 
\mathtt{Normalize}(\mathcal{A}^T*\mathcal{M}^{-1}*\mathcal{\vec{Q}}_1, \mathcal{L})$ 
with $\mathbf{c}_1$ invertible\\
\For {$i = 1,2, \dots, k$}{
$\mathcal{\vec{Q}} \leftarrow \mathcal{A*L}*\mathcal{\vec{W}}_i-
\mathcal{\vec{Q}}_i*\mathbf{c}_i $\\
$[\mathcal{\vec{Q}}_{i+1}, {\mathbf{z}}_{i+1}] \leftarrow 
\mathtt{Normalize}(\mathcal{\vec{Q}}, \mathcal{M}^{-1})$\\
\If {$i<k$}{
$\mathcal{\vec{W}} \leftarrow \mathcal{A}^T*\mathcal{M}^{-1}*\mathcal{\vec{Q}}_{i+1}-
\mathcal{\vec{W}}_i* \mathbf{z}_{i+1}$\\
$[\mathcal{\vec{W}}_{i+1}, \mathbf{c}_{i+1}] \leftarrow 
\mathtt{Normalize}(\mathcal{\vec{W}}, \mathcal{L})$\\
}
}
\caption{The partial weighted tensor Golub-Kahan bidiagonalization (W-tGKB) process}\label{Alg: 3}
\end{algorithm}\vspace{.3cm}

Algorithm \ref{Alg: 3} produces the partial W-tGKB decompositions
\begin{equation}\label{4.1}
\mathcal{A*L} * \mathcal{W}_k = \mathcal{Q}_{k+1} * \mathcal{\bar{P}}_k, \;\;\;\;\;
\mathcal{A}^T *\mathcal{M}^{-1}* \mathcal{Q}_k = \mathcal{W}_k * \mathcal{P}_k^T,
\end{equation}
where
\[
\mathcal{\bar{P}}_k = \begin{bmatrix}
\mathbf{c}_1 & & &\\
\mathbf{z}_2 & \mathbf{c}_2 & \\
& \mathbf{z}_3 & \mathbf{c}_3 & \\
& & \ddots & \ddots & \\
&	& & \mathbf{z}_{k} & \mathbf{c}_k\\
&	& & & \mathbf{z}_{k+1} 
\end{bmatrix} \in \mathbb{R}^{(k+1)\times k\times n}
\]
is a lower bidiagonal tensor and $\mathcal{P}_k$ denotes the leading $k\times k\times n$ 
subtensor of $\mathcal{\bar{P}}_k$. The tensor columns 
$\mathcal{\vec{W}}_j\in\mathbb{R}^{m\times 1\times n}$, $j=1,2,\dots,k$, and
$\mathcal{\vec{Q}}_{i+1}\in\mathbb{R}^{\ell\times 1\times n}$, $i=0,1,\dots,k$,
generated by Algorithm \ref{Alg: 3} are orthonormal tensor bases for the t-Krylov 
subspaces 
\[
\mathbb{K}_k(\mathcal{A}^T*\mathcal{M}^{-1}*\mathcal{A*L},
\mathcal{A}^T*\mathcal{M}^{-1}*\mathcal{\vec{B}})~~~\mbox{and}~~~ 
\mathbb{K}_{k+1}(\mathcal{A*L}*\mathcal{A}^T*\mathcal{M}^{-1},\mathcal{\vec{B}}),
\]
respectively. The tensor columns $\mathcal{\vec{W}}_j$ and $\mathcal{\vec{Q}}_i$ define 
the tensors
\begin{equation*}
\mathcal{W}_k := [\mathcal{\vec{W}}_1, \dots, \mathcal{\vec{W}}_k] \in 
\mathbb{R}^{m \times k \times n}, ~~~\mathcal{Q}_{k+1} := [\mathcal{\vec{Q}}_1, \dots, 
\mathcal{\vec{Q}}_{k+1}] \in \mathbb{R}^{\ell \times (k+1) \times n},
\end{equation*}
which satisfy 
\be \label{normLM}
\mathcal{W}_k^T*\mathcal{L}*\mathcal{W}_k = \mathcal{I}_k ~~{\rm and}~~ 
\mathcal{Q}_{k+1}^T*\mathcal{M}^{-1}*\mathcal{Q}_{k+1} = \mathcal{I}_{k+1}.
\ee

Theorem \ref{spatial} below shows that the application of the operators $\mathcal{L}$ and 
$\mathcal{M}$ in the spatial domain is equivalent to their application in the Fourier 
domain. Application of these operators in the spatial domain may reduce the computing 
time since transformation of $\mathcal{L}$ and $\mathcal{M}$ to and from the Fourier 
domain is not required. We apply the SPD tensors $\mathcal{L}$ and $\mathcal{M}$ in the 
spatial domain in the computed examples of Section \ref{sec5}.

\begin{thm}\label{spatial}
Let the first frontal slice $\mathcal{L}^{(1)}$ of the tensor 
$\mathcal{L}\in\mathbb{R}^{s\times m\times n}$ be nonzero and let the remaining frontal 
slices, $\mathcal{L}^{(i)}$, $i=2,3,\dots,n$, be zero matrices. Let 
$\mathcal{X}\in\mathbb{R}^{m\times p\times n}$ have frontal slices $\mathcal{X}^{(i)}$, 
$i=1,2,\dots,n$. Define
\[
\mathcal{L}(\mathcal{X}) := \mathcal{L*X}.
\]
Then
\[
\mathcal{L}(\mathcal{X})^{(i)} = \mathcal{L}^{(1)}\mathcal{X}^{(i)}, ~~~ i=1,2,\dots,n.
\]
\end{thm}
\noindent
{\it Proof:} We have
\begin{equation}\label{spatialw}
\begin{split}
\mathcal{L*X} &= {\tt fold}\left({\tt bcirc}(\mathcal{L}) 
\cdot{\tt unfold}(\mathcal{X})\right) \\ &= {\tt fold}\left( \begin{bmatrix}
\mathcal{L}^{(1)} & O & \dots & O\\
O &\mathcal{L}^{(1)} & \dots & O\\
\vdots & \vdots & \ddots & \vdots \\
O & O & \dots & \mathcal{L}^{(1)}
\end{bmatrix} \begin{bmatrix}
\mathcal{X}^{(1)}\\
\mathcal{X}^{(2)}\\
\vdots\\
\mathcal{X}^{(n)}
\end{bmatrix}\right) \\ &= {\tt fold}\left(\begin{bmatrix}
\mathcal{L}^{(1)} \mathcal{X}^{(1)}\\
\mathcal{L}^{(1)} \mathcal{X}^{(2)}\\
\vdots\\
\mathcal{L}^{(1)} \mathcal{X}^{(n)}
\end{bmatrix} \right),
\end{split}
\end{equation}
and the proof follows. ~~~$\Box$

\subsection{The W-tGKT method}\label{sec3.1}
This subsection describes the W-tGKT method for the approximate solution of least squares 
problems of the form \eqref{4w}. We also apply this method $p$ times to determine an 
approximate solution of \eqref{4m} by working with the data tensor slices 
$\mathcal{\vec{B}}_j$, $j=1,2,\dots,p$, of $\mathcal{B}$ independently. 

Let $\mathcal{\vec{Y}} = \mathcal{W}_k*\mathcal{\vec{Z}}$. Substituting the left-hand 
W-tGKB decomposition \eqref{4.1} into \eqref{n4} and using \eqref{normLM} yields
\begin{equation}
\min_{\mathcal{\vec{Z}} \in \mathbb{R}^{k \times 1 \times n}}
\{ \|\mathcal{\bar{P}}_{k}*\mathcal{\vec{Z}} - \mathcal{Q}_{k+1}^T*\mathcal{\vec{B}}\|^2_F
+ \mu^{-1}\|\mathcal{\vec{Z}}\|^2_F\}.
\label{4.2}
\end{equation}
Since $\mathcal{\vec{B}} = \mathcal{\vec{Q}}_1*\mathbf{z}_1$ (cf. Algorithm \ref{Alg: 3}),
we get 
\begin{equation}\label{QTB}
\mathcal{Q}_{k+1}^T*\mathcal{\vec{B}}=\vec{e}_1*\mathbf{z}_1\in
\mathbb{R}^{(k+1)\times 1\times n},
\end{equation}
where the $(1,1,1)$th entry of $\vec{\mathit{e}}_1\in\mathbb{R}^{m\times 1\times n}$ is 
$1$ and the remaining entries are zero. Substituting \eqref{QTB} into \eqref{4.2}, we 
obtain
\begin{equation}
\min_{\mathcal{\vec{Z}}\in \mathbb{R}^{k \times 1 \times n}}
\{ \|\mathcal{\bar{P}}_{k}*\mathcal{\vec{Z}} - \vec{\mathit{e}}_1*\mathbf{z}_1\|^2_F +
\mu^{-1}\|\mathcal{\vec{Z}}\|^2_F\}.
\label{4.4}
\end{equation}
The equivalency of \eqref{4w} and \eqref{4.4} follows from \eqref{normLM}. Thus,
\begin{equation} \label{4.4equ}
\min_{\mathcal{\vec{Y}}\in \mathbb{K}_k}
\{\|\mathcal{A*L*\vec{Y}}-\mathcal{\vec{B}}\|^2_{\mathcal{M}^{-1}} + 
\mu^{-1}\|\mathcal{\vec{Y}}\|^2_{\mathcal{L}}\} \iff 
\min_{\mathcal{\vec{Z}}\in \mathbb{R}^{k \times 1 \times n}}
\{ \|\mathcal{\bar{P}}_{k}*\mathcal{\vec{Z}} - \vec{\mathit{e}}_1*\mathbf{z}_1\|^2_F +
\mu^{-1}\|\mathcal{\vec{Z}}\|^2_F\},
\end{equation}
where $\mathbb{K}_k=\mathbb{K}_k(\mathcal{A}^T*\mathcal{M}^{-1}*\mathcal{A*L},
\mathcal{A}^T*\mathcal{M}^{-1}*\mathcal{\vec{B}})$.

A suitable way to compute the solution of the minimization problem on the right-hand side 
of \eqref{4.4equ} is to solve the least squares problem
\begin{equation}\label{4.9}
\min_{\mathcal{\vec{Z}}\in \mathbb{R}^{k \times 1 \times n}}
\left\|\left[\begin{array}{c} \mathcal{\bar{P}}_k \\
\mu^{-1/2}\mathcal{I} \end{array} \right] * \mathcal{\vec{Z}} -
\left[\begin{array}{c}
\vec{\mathit{e}}_1*\mathbf{z}_1\\
\mathcal{\vec{O}}
\end{array}\right]\right\|_F,
\end{equation}
using \cite[Algorithm 3]{RU1}. The solution $\mathcal{\vec{Z}}_{\mu,k}$ can be expressed 
as
\begin{equation}\label{4.10}
\mathcal{\vec{Z}}_{\mu,k} = (\mathcal{\bar{P}}_k^T*\mathcal{\bar{P}}_k +
\mu^{-1} \mathcal{I})^{-1}*\mathcal{\bar{P}}_k^T*{\vec{\mathit{e}}}_1*\mathbf{z}_1,
\end{equation}
and the associated approximate solution of \eqref{4w} is given by
\[
\mathcal{\vec{X}}_{\mu,k} = \mathcal{L}*\mathcal{W}_k *\mathcal{\vec{Z}}_{\mu,k}.
\]

The discrepancy principle \eqref{discr} is used to determine the regularization parameter
and the number of steps of the W-tGKB process. We show below that \eqref{discr} is 
equivalent to 
\begin{equation} \label{41.22n}
\|\mathcal{\bar{P}}_{k}*\mathcal{\vec{Z}}_{\mu,k}-
\vec{\mathit{e}}_1*\mathbf{z}_1\|_F = \eta \delta.
\end{equation}
Thus, we have to determine $\mu>0$ so that \eqref{41.22n} holds. Define the function
\be \label{41.22}
\phi_k(\mu):=\|\mathcal{\bar{P}}_{k}*\mathcal{\vec{Z}}_{\mu,k}-
\vec{\mathit{e}}_1*\mathbf{z}_1\|_F^2.
\ee
Substituting \eqref{4.10} into \eqref{41.22}, and using
the identity
\be \label{iden}
\mathcal{I}-\mathcal{\bar{P}}_k*(\mathcal{\bar{P}}_k^T*\mathcal{\bar{P}}_k+\mu^{-1}
\mathcal{I})^{-1}*\mathcal{\bar{P}}_k^T = (\mu\mathcal{\bar{P}}_k*\mathcal{\bar{P}}_k^T +\mathcal{I})^{-1}
\ee
as well as \eqref{defn: 2.7}, we obtain
\[
\phi_k(\mu) = \big((\vec{\mathit{e}}_1*\mathbf{z}_1)^T*(\mu\mathcal{\widetilde{P}}_k*
\mathcal{\widetilde{P}}_k^T+\mathcal{I})^{-2}*\vec{\mathit{e}}_1*
\mathbf{z}_1\big)_{(:,:,1)}.
\]
Therefore, \eqref{41.22n} becomes
\be \label{4.14}
\phi_k(\mu) - (\eta \delta)^2 = 0.
\ee

The dependence of $\phi_k(\mu)$ on $k$ for a fixed $\mu>0$ is established in \cite{RU1}.
Moreover, $\phi_k(\mu)$ is shown to be a decreasing and convex function of $\mu>0$. We can
define $\phi_k(0) = \|\vec{\mathit{e}}_1*\mathbf{z}_1\|_F^2$ by continuity. In the 
computed examples of Section \ref{sec5}, we use the bisection method to solve 
\eqref{4.14}. The following result establishes the equivalence between the discrepancy
principle \eqref{discr} and \eqref{41.22n}.

\begin{prop}\label{prop32}
Let $\phi_k(\mu)$ be defined by \eqref{41.22}, and assume that $\mu=\mu_k$ solves
$\phi_k(\mu)=\eta^2\delta^2$ and that $\mathcal{\vec{Z}}_{\mu,k}$ solves
\begin{equation*}
(\mathcal{\bar{P}}_k^T*\mathcal{\bar{P}}_k+\mu^{-1}\mathcal{I})*
\mathcal{\vec{Z}}=\mathcal{\bar{P}}_k^T*{\vec{\mathit{e}}}_1*\mathbf{z}_1.
\end{equation*}
Let $\mathcal{\vec{Y}}_{\mu,k} =\mathcal{W}_k*\mathcal{\vec{Z}}_{\mu,k}$.
Then the associated approximate solution
$\mathcal{\vec{X}}_{\mu,k} = \mathcal{L}*\mathcal{\vec{Y}}_{\mu,k}$ of \eqref{srhs} satisfies
\[
\|\mathcal{A}*\mathcal{\vec{X}}_{\mu,k} - \mathcal{\vec{B}}\|_{\mathcal{M}^{-1}}^2=
\big((\vec{\mathit{e}}_1*\mathbf{z}_1)^T*(\mu\mathcal{\bar{P}}_k*
\mathcal{\bar{P}}_k^T + \mathcal{I})^{-2}*\vec{\mathit{e}}_1*
\mathbf{z}_1\big)_{(:,:,1)}.
\]
\end{prop}

\noindent
\textit{Proof}: Substituting
$\mathcal{\vec{X}}_{\mu,k}=\mathcal{L}*\mathcal{\vec{Y}}_{\mu,k}$ and $\mathcal{\vec{Y}}_{\mu,k} =\mathcal{W}_k*\mathcal{\vec{Z}}_{\mu,k}$ into \eqref{discr},
and using the left-hand side decomposition of \eqref{4.1}, as well as \eqref{QTB} and \eqref{normLM}, shows
that
\begin{equation*}
\|\mathcal{A}*\mathcal{\vec{X}}_{\mu,k}-\mathcal{B}\|_{\mathcal{M}^{-1}}^2=
\|\mathcal{Q}_{k+1}*\mathcal{\bar{P}}_k*\mathcal{\vec{Z}}_{\mu,k}-
\mathcal{\vec{B}}\|_{\mathcal{M}^{-1}}^2=
\|\mathcal{\bar{P}}_k*\mathcal{\vec{Z}}_{\mu,k}-\vec{e}_1*\mathbf{z}_1\|_F. ~~~\Box
\end{equation*}

We refer to the solution method for \eqref{4w} described above as the W-tGKT method.
It is implemented by Algorithm \ref{Alg: 6} below with $p=1$.

The remainder of this subsection describes an algorithm for the approximate solution of 
\eqref{4m} by solving the minimization problem \eqref{nn4}. Let the data tensor 
$\mathcal{B}$ have the lateral data slices $\mathcal{\vec{B}}_j$, $j=1,2,\ldots,p$, and 
let $\mathcal{\vec{B}}_{j,{\rm true}}$ denote the unknown error-free data slice 
associated with the available error-contaminated slice $\mathcal{\vec{B}}_j$. Assume that 
bounds for the norm of the errors
\[
\mathcal{E}_j:=\mathcal{\vec{B}}_j-\mathcal{\vec{B}}_{j,{\rm true}},\quad j=1,2,\dots,p, 
~~~p>1,
\]
are available or can be estimated, i.e.,
\[
\|\mathcal{\vec{E}}_j\|_{\mathcal{M}^{-1}}\leq\delta_j,\quad j=1,2,\ldots,p;
\]
cf. \eqref{Btrue} and \eqref{errbd}. Let 
$\mathcal{\vec{X}}_j = \mathcal{L}*\mathcal{\vec{Y}}_j$, $j = 1,2,\dots,p$, and consider 
\eqref{nn4} as $p$ separate Tikhonov minimization problems
\begin{equation}\label{33n}
\min_{\mathcal{\vec{Y}}_j \in \mathbb{R}^{m\times 1\times n}}
\{\|\mathcal{A*L}*\mathcal{\vec{Y}}_j-\mathcal{\vec{B}}_j\|_{\mathcal{M}^{-1}}^2+
{\mu}^{-1}\|\mathcal{\vec{Y}}_j\|_\mathcal{L}^2\}, \;\;\;\; j = 1,2,\dots,p.
\end{equation}
We apply the W-tGKT method described above to solve the $p$ problems \eqref{33n} 
independently by using Algorithm \ref{Alg: 6}, and we refer to this solution method for 
\eqref{4m} by solving the problems \eqref{33n} independently as the W-tGKT$_p$ method.

\vspace{.3cm}
\begin{algorithm}[H]
\SetAlgoLined
\KwIn{$\mathcal{A}, \mathcal{\vec{B}}_1, \mathcal{\vec{B}}_2, \dots, \mathcal{\vec{B}}_p,
\delta_1, \delta_2, \dots, \delta_p, \mathcal{L}, \mathcal{M}, \eta > 1, k_\text{init}= 2$}
\For{$j = 1,2,\dots, p$}{
$k \leftarrow k_\text{init}$, $[\mathcal{\vec{Q}}_1, \mathbf{z}_1] \leftarrow \mathtt{Normalize}(\mathcal{\vec{B}}_j,\mathcal{M}^{-1})$\\
Compute $ \mathcal{W}_k, \mathcal{Q}_{k+1}$ and $\mathcal{\mathcal{\bar{P}}}_k$ by Algorithm \ref{Alg: 3}, and let $\vec{\mathit{e}}_1 \leftarrow \mathcal{I}(:,1,:)$\\
Solve the minimization problem
\[
\min_{\mathcal{\vec{Z}}\in \mathbb{R}^{k \times 1 \times n}}
\|\mathcal{\bar{P}}_k*\mathcal{\vec{Z}} - \vec{\mathit{e}}_1*\mathbf{z}_1 \|_F
\]
for $\mathcal{\vec{Z}}_k$ by using \cite[Algorithm 3]{RU1}\\

\While{$\|\mathcal{\bar{P}}_k*\mathcal{\vec{Z}}_k - \vec{\mathit{e}}_1*\mathbf{z}_1\|_F \geq \eta \delta_j$}{
$k \leftarrow k+1$\\
$\mathtt{Go \; to \; step \; 3}$}
Determine the regularization parameter by the discrepancy principle, i.e., compute the zero $\mu_k>0$ of
\[
\xi_k(\mu_k):=\|\mathcal{\bar{P}}_k*\mathcal{\vec{Z}}_{j,\mu_k} -
\vec{\mathit{e}}_1*\mathbf{z}_1\|_F^2-\eta^2 \delta_j^2
\]
and the associated solution $\mathcal{\vec{Z}}_{j,\mu_k}$ of
\[
\min_{\mathcal{\vec{Z}}\in \mathbb{R}^{k \times 1 \times n}}\left\|\left[\begin{array}{c}
\mathcal{\bar{P}}_k \\
\mu_k^{-1/2}\mathcal{I}
\end{array}\right] * \mathcal{\vec{Z}} - \left[\begin{array}{c}
\vec{\mathit{e}}_1*\mathbf{z}_1\\
\mathcal{\vec{O}}
\end{array}\right] \right\|_F
\]
by using \cite[Algorithm 3]{RU1} \\
Compute $\mathcal{\vec{Y}}_{j,\mu_k} \leftarrow \mathcal{W}_k*\mathcal{\vec{Z}}_{j,\mu_k},
\;\; \mathcal{\vec{X}}_{j,\mu_k} \leftarrow \mathcal{L}*\mathcal{\vec{Y}}_{j,\mu_k}$\\
}
\caption{The W-tGKT$_p$ method for the approximate solution of \eqref{4m} by solving the
problems \eqref{33n} independently}\label{Alg: 6}
\end{algorithm} \vspace{.3cm}

\section{The WGG-tGKT and WG-tGKT Methods}\label{sec4}
We describe the WGG-tGKT and WG-tGKT methods for the approximate solution of \eqref{4w} 
and \eqref{4m}. The latter method is designed to compute an approximate solution of 
\eqref{4w}. When applied to solve \eqref{4m}, the method works with the lateral slices of 
$\mathcal{\vec{B}}_j$, $j=1,2,\dots,p$, of $\mathcal{B}$ one at a time and independently. 
We refer to the resulting method for solving \eqref{4m} as the WG-tGKT$_p$ method. The 
WGG-tGKT method solves \eqref{4m} by working with the lateral slices of the tensor 
$\mathcal{B}$ simultaneously.

\subsection{The WGG-tGKT method}\label{sec4.1}
This subsection describes the weighted generalized global t-product Golub-Kahan 
bidiagonalization (WGG-tGKB) process and discusses how it can be applied to reduce the 
large-scale problem \eqref{4m} to a problem of small size. This is achieved by applying 
$k \ll m$ steps of the WGG-tGKB process to $\mathcal{A}$. We assume that there is no 
breakdown of this process. This is the generic situation. The WGG-tGKB process, which is
described by Algorithm \ref{Alg: 8}, yields the decompositions
\begin{equation}\label{decom}
\mathcal{A}*\mathcal{L}*\mathbb{W}_k = \mathbb{Q}_{k+1} \circledast \bar{P}_k, \quad 
\mathcal{A}^T*\mathcal{M}^{-1}*\mathbb{Q}_k = \mathbb{W}_k \circledast P_k^T,
\end{equation}
where
\begin{equation*}
\mathbb{W}_k := [\mathcal{W}_1,\dots,\mathcal{W}_k]\in\mathbb{R}^{m \times kp \times n},
\quad \mathbb{Q}_{k+1} := [\mathcal{Q}_1,\dots,\mathcal{Q}_{k+1}]\in
\mathbb{R}^{\ell \times (k+1)p \times n}.
\end{equation*}
As usual, the tensors $\mathcal{L}\in\mathbb{R}^{m \times m\times n}$ and 
$\mathcal{M}\in\mathbb{R}^{\ell\times\ell\times n}$ are assumed to be symmetric positive 
definite. Moreover,
\be \label{Snorm}
\mathbb{Q}_{k+1}^T\Diamond_{\mathcal{M}^{-1}} \mathbb{Q}_{k+1} = I_{k+1},\quad
\mathbb{W}_k^T\Diamond_\mathcal{L}\mathbb{W}_k = I_k,
\ee
and
\begin{equation}\label{spec}
\begin{split}
\mathcal{A*L}*\mathbb{W}_k&=[\mathcal{A*L}*\mathcal{W}_1,\mathcal{A*L}*\mathcal{W}_2,\dots,
\mathcal{A*L}*\mathcal{W}_k]\in\mathbb{R}^{\ell \times kp \times n},\\
\mathbb{Q}_{k+1}\circledast\bar{P}_k&=[\mathbb{Q}_{k+1}\circledast \bar{P}_k(:,1),
\mathbb{Q}_{k+1}\circledast \bar{P}_k(:,2),\dots,\mathbb{Q}_{k+1}\circledast
\bar{P}_k(:,k)]\in\mathbb{R}^{\ell\times kp\times n},
\end{split}
\end{equation}
where $\mathcal{A}^T*\mathcal{M}^{-1}*\mathbb{Q}_k$ and 
$\mathbb{W}_k \circledast P_k^T$ are defined similarly as \eqref{spec}.

Details of the computations of the WGG-tGKB process are described by Algorithm 
\ref{Alg: 8}. The tensors $\mathcal{W}_j\in\mathbb{R}^{m\times p\times n}$, 
$j=1,2,\dots,k$, and $\mathcal{Q}_{i+1}\in\mathbb{R}^{\ell\times p\times n}$, 
$i=0,1,\dots,k$, generated by the algorithm form orthogonal tensor bases for the t-Krylov
subspaces
$\mathbb{K}_k(\mathcal{A}^T*\mathcal{M}^{-1}*\mathcal{A*L}, 
\mathcal{A}^T*\mathcal{M}^{-1}*\mathcal{B})$ and 
$\mathbb{K}_{k+1}(\mathcal{A*L}*\mathcal{A}^T*\mathcal{M}^{-1},\mathcal{B})$, 
respectively. The lower bidiagonal matrix $\bar{P}_k$ in \eqref{spec} is given by
\begin{equation}\label{tridiag}
\bar{P}_k = \begin{bmatrix}
\alpha_1 & & &\\
\beta_2 & \alpha_2 & \\
& \beta_3 & \alpha_3 & \\
& & \ddots & \ddots & \\
& & & \beta_{k} & \alpha_k\\
& & & & \beta_{k+1}
\end{bmatrix}\in\mathbb{R}^{(k+1)\times (k+1)},
\end{equation}
and $P_k$ is the leading $k\times k$ submatrix of $\bar{P}_k$. The relation
\be
\mathcal{B} = \mathbb{Q}_{k+1} \circledast e_1 \beta_1, \;\;\; e_1 = [1,0,\dots,0]^T,
\label{A352}
\ee
is easily deduced from Algorithm \ref{Alg: 8}. 

\begin{prop}\label{prop3m}
Let the tensors $\mathbb{Q}_{k+1}$ and $\mathcal{M}$ be defined as above. Let 
$y\in\mathbb{R}^{k+1}$. Then
\be
\|\mathbb{Q}_{k+1}\circledast y\|_{\mathcal{M}^{-1}}=\|y\|_2,
\ee
where $\|\cdot \|_2$ denotes the Euclidean vector norm.
\end{prop}
\noindent
{\it Proof:}
\[
\|\mathbb{Q}_{k+1}\circledast y\|_{\mathcal{M}^{-1}}^2 = 
\left \langle \sum_{i=1}^{k+1} y_j \mathcal{Q}_j, 
\sum_{i=1}^{k+1} y_j \mathcal{Q}_j \right \rangle_{\mathcal{M}^{-1} } = 
\sum_{i=1}^{k+1} y_j^2 \langle \mathcal{Q}_j, \mathcal{Q}_j \rangle_{\mathcal{M}^{-1} } = 
\sum_{i=1}^{k+1} y_j^2 = \|y\|_2^2,
\]
since the tensors $\mathcal{Q}_j$, $j=1,2,\dots,k$, are orthogonal, i.e.,
\[
\langle \mathcal{Q}_i, \mathcal{Q}_j \rangle_{\mathcal{M}^{-1} } = \left\{
\begin{array}{ll}
1 & i=j,\\
0 & i \neq j.
\end{array}
\right. 
\]
\hfill$\Box$\break

It can be shown analogously that for an orthogonal tensor 
$\mathcal{Q}\in\mathbb{R}^{m\times k\times n}$, one has
\begin{equation}\label{norm2F}
\|\mathcal{Q}\circledast y\|_{\mathcal{M}^{-1}}=\|y\|_2;
\end{equation}
see \cite{GIJS} for the analogous result when $\mathcal{M}^{-1} = \mathcal{I}$.

We compute an approximate solution of \eqref{4m} similarly as described in Subsection 
\ref{sec3.1}. Thus, letting $\mathcal{Y} = \mathbb{W}_k \circledast z$ and using the 
left-hand side of \eqref{decom}, as well as \eqref{A352} and \eqref{Snorm}, the 
minimization problem \eqref{nn4} reduces to
\begin{equation}\label{N1}
\min_{z\in\mathbb{R}^k}\{\|\bar{P}_k z- e_1\beta_1\|^2_2+\mu^{-1}\|z\|^2_2\}.
\end{equation}

\vspace{.3cm}
\begin{algorithm}[H]
\SetAlgoLined
\KwIn{$ \mathcal{A} \in \mathbb{R}^{\ell \times m \times n}$, $\mathcal{B} \in
\mathbb{R}^{\ell \times p \times n}$ such that 
$\mathcal{A}^T*\mathcal{M}^{-1}*\mathcal{B}\neq\mathcal{O}$, 
$\mathcal{L} \in \mathbb{R}^{m \times m \times n}$, $\mathcal{M} \in 
\mathbb{R}^{\ell \times \ell \times n}$ are SPD, $k\geq 1$}
$\beta_1 \leftarrow \|\mathcal{B}\|_{\mathcal{M}^{-1}}$, $\mathcal{Q}_1 \leftarrow
\frac{1}{\beta_1}\mathcal{B}$\\
$\alpha_1 \leftarrow \|\mathcal{A}^T*\mathcal{M}^{-1}*\mathcal{Q}_1\|_\mathcal{L}$, $\mathcal{W}_1 \leftarrow \frac{1}{\alpha_1}(\mathcal{A}^T*\mathcal{M}^{-1}*\mathcal{Q}_1)$\\
\For {$j = 1,2, \dots, k$}{
$\mathcal{Q} \leftarrow \mathcal{A*L}*\mathcal{W}_j - \alpha_j\mathcal{Q}_j$\\
$\beta_{j+1} \leftarrow \|\mathcal{Q}\|_{\mathcal{M}^{-1}}$, 
$\mathcal{Q}_{j+1} \leftarrow \mathcal{Q}/\beta_{j+1}$\\
\If {$j<k$}{
$\mathcal{W} \leftarrow \mathcal{A}^T*\mathcal{M}^{-1}*\mathcal{Q}_{j+1} - \beta_{j+1}\mathcal{W}_j$\\
$\alpha_{j+1} \leftarrow \|\mathcal{W}\|_\mathcal{L}$, $\mathcal{W}_{j+1} \leftarrow \mathcal{W}/\alpha_{j+1}$\\
}
}
\caption{Partial weighted generalized global tGKB process}
\label{Alg: 8}
\end{algorithm}\vspace{.3cm}

The minimization problem \eqref{N1} is analogous to \eqref{4.9}. We compute its solution 
by solving
\begin{equation}
\min_{z\in \mathbb{R}^k}\left\| \left[\begin{array}{c}
\bar{P}_k \\
\mu^{-1/2} I
\end{array}\right] z - \left[\begin{array}{c}
e_1 \beta_1\\
0
\end{array}\right] \right\|_2.
\label{A6}
\end{equation}
Denote the solution by $z_{\mu,k}$. Then the associated approximate solution of
\eqref{4m} is given by
\[
\mathcal{\vec{X}}_{\mu,k} = \mathcal{L}*\mathbb{W}_k \circledast z_{\mu,k}.
\]

We determine the regularization parameter $\mu$ by the discrepancy principle based on the
$\mathcal{M}^{-1}$-norm. This assumes knowledge of a bound
\[
\|\mathcal{E}\|_{\mathcal{M}^{-1}}\leq\delta
\]
for the error tensor $\mathcal{E}$ in $\mathcal{B}$. Thus, we choose $\mu>0$ so that the 
solution $z_{\mu,k}$ of \eqref{A6} satisfies
\begin{equation}
\|\bar{P}_k z_{\mu,k} - e_1 \beta_1\|_2 = \eta \delta.
\label{AA9}
\end{equation}

Define the function
\[
\psi_k(\mu):=\|\bar{P}_kz_{\mu,k}-e_1 \beta_1\|_2^2,
\]
where $z_{\mu,k}$ solves \eqref{A6}. Using an identity analogous to \eqref{iden}, we 
obtain
\[
\psi_k(\mu) = \beta_1^2e_1^T(\mu \bar{P}_k\bar{P}_k^T + I)^{-2}e_1,
\]
The discrepancy principle \eqref{AA9} can be satisfied for reasonable values of 
$\eta\delta$ by choosing a large enough value of $k$; see \cite[Proposition 4.3]{RU1}. It 
can be shown that the function $\mu\mapsto\psi_k(\mu)$ is decreasing and convex with 
$\psi_k(0)=\beta_1^2$. A matrix analogue of Proposition \ref{prop32} is shown below.

\begin{prop}\label{prop3d}
Let $\mu_k$ solve $\psi_k(\mu)=\eta^2\delta^2$ and suppose that $z_{\mu,k}$ is the
solution of \eqref{N1} with $\mu=\mu_k$. Let $\mathcal{Y}_{\mu,k} = \mathbb{W}_k \circledast z_{\mu,k}$. Then the associated approximate solution
$\mathcal{X}_{\mu,k} = \mathcal{L} * \mathcal{Y}_{\mu,k}$ of \eqref{4m} satisfies
\begin{equation}
\|\mathcal{A}*\mathcal{X}_{\mu,k} - \mathcal{B}\|_{\mathcal{M}^{-1}}^2 =
\beta_1^2e_1^T(\mu \bar{P}_k\bar{P}_k^T + I)^{-2}e_1.
\label{4.12}
\end{equation}
\end{prop}
\noindent
{\it Proof:} Substituting $\mathcal{X}_{\mu,k}= \mathcal{L} * \mathcal{Y}_{\mu,k}$ and $ \mathcal{Y}_{\mu,k} = \mathbb{W}_k \circledast z_{\mu,k}$ into left-hand side
of \eqref{4.12}, using the left-hand side of \eqref{decom}, \eqref{Snorm} and \eqref{A352}, as well as left-hand side of
\eqref{norm2F} gives
\begin{equation*}
\|\mathcal{A}*\mathcal{X}_{\mu,k} - \mathcal{B}\|_{\mathcal{M}^{-1}}^2
= \|\mathbb{Q}_{k+1} \circledast (\bar{P}_k \circledast z_{\mu,k} - e_1 \beta_1)\|_{\mathcal{M}^{-1}}^2
= \| \bar{P}_k z_{\mu,k} - e_1 \beta_1\|_2^2. \;\; \Box
\end{equation*}
We refer to the solution method described above as the WGG-tGKT method. This method is
implemented by Algorithm \ref{Alg: 10}.

\vspace{.3cm}
\begin{algorithm}[H]
\SetAlgoLined
\KwIn{$\mathcal{A}$, $\mathcal{B}$, $\delta$, $\mathcal{L}$, $\mathcal{M}$, $\eta > 1$, $ k_\text{init} = 2$}
$k \leftarrow k_\text{init}$, $\beta_1 \leftarrow \|\mathcal{B}\|_{\mathcal{M}^{-1}}$,
$\mathcal{Q}_1 \leftarrow \frac{1}{\beta_1}\mathcal{B}$\\
Compute $\mathbb{W}_k$, $\mathbb{Q}_{k+1}$, and $\bar{P}_k$ by Algorithm \ref{Alg: 8}\\
Let $z_k\in\mathbb{R}^k$ solve the minimization problem
\[
\min_{z \in \mathbb{R}^k} \| \bar{P}_k z - e_1 \beta_1 \|_2
\]
\While{$\|\bar{P}_k z_k - e_1 \beta_1\|_2 \geq \eta \delta$}{
$k \leftarrow k+1$\\
$\mathtt{Go \; to \; step}$ 2
}
Determine the regularization parameter by the discrepancy principle, i.e., compute the zero $\mu_k>0$ of
\[
\varphi_k(\mu):=\|\bar{P}_k z_{\mu,k} - e_1 \beta_1\|_2^2-\eta^2 \delta^2
\]
and the associated solution $z_{\mu,k}$ of
\[
\min_{z \in \mathbb{R}^k}\left\| \left[\begin{array}{c}
\bar{P}_k \\
\mu_k^{-1/2} I
\end{array}\right] z - \left[\begin{array}{c}
e_1 \beta_1\\
0
\end{array}\right] \right\|_2
\] \\
Compute $ \mathcal{Y}_{\mu,k} \leftarrow \mathbb{W}_k\circledast z_{\mu,k}, \;\; \mathcal{X}_{\mu,k} \leftarrow \mathcal{L}*\mathcal{Y}_{\mu,k}$ \\
\caption{The WGG-tGKT method for computing an approximate solution of \eqref{4m}}
\label{Alg: 10}
\end{algorithm} \vspace{.3cm}

\subsection{The WG-tGKT method}\label{sec4.2}
The WG-tGKT method for the approximate solution of \eqref{4w} first reduces the tensor 
$\mathcal{A}$ in \eqref{4w} to a small lower bidiagonal matrix by carrying out a few, say
$k$, steps of the weighted global t-product Golub-Kahan bidiagonalization (WG-tGKB)
process, which is described by Algorithm \ref{Alg: 11}. We assume that $k$ is small enough
to avoid breakdown. This is the generic situation. Algorithm \ref{Alg: 11} yields the 
partial WG-tGKB decompositions
\be \label{BB}
\mathcal{A*L}*\mathcal{W}_k = \mathcal{Q}_{k+1} \circledast \bar{B}_k, \quad
\mathcal{A}^T*\mathcal{M}^{-1}*\mathcal{Q}_k = \mathcal{W}_k \circledast B_k^T,
\ee
where
\[
\mathcal{W}_k := [\mathcal{\vec{W}}_1,\dots,\mathcal{\vec{W}}_k] \in
\mathbb{R}^{m \times k \times n}, \;\;
\mathcal{Q}_{k+1} := [\mathcal{\vec{Q}}_1, \dots, \mathcal{\vec{Q}}_{k+1}]\in
\mathbb{R}^{\ell \times (k+1) \times n}
\]
and
\[
\mathcal{Q}_{k+1}^T\Diamond_{\mathcal{M}^{-1}} \mathcal{Q}_{k+1} = I_{k+1}, \quad
\mathcal{W}_k^T\Diamond_\mathcal{L}\mathcal{W}_k = I_k.
\]
The tensors $\mathcal{L}$ and $\mathcal{M}$ are assumed to be SPD. The expressions in 
\eqref{BB} are defined similarly to \eqref{spec}, and the bidiagonal matrix 
$\bar{B}_k \in \mathbb{R}^{(k+1) \times k}$ has a form analogous to \eqref{tridiag}. The 
tensors $\mathcal{\vec{W}}_j\in\mathbb{R}^{\ell\times 1\times n}$, $j = 1,2,\dots,k$, 
and $\mathcal{\vec{Q}}_{i+1}\in\mathbb{R}^{m\times 1\times n}$, $i=0,1,\dots,k$, generated 
by Algorithm \ref{Alg: 11} form orthonormal tensor bases for the t-Krylov subspaces
$\mathbb{K}_k(\mathcal{A}^T*\mathcal{M}^{-1}*\mathcal{A*L},
\mathcal{A}^T*\mathcal{M}^{-1}*\mathcal{\vec{B}})$ and 
$\mathbb{K}_{k+1}(\mathcal{A*L}*\mathcal{A}^T*\mathcal{M}^{-1},\mathcal{\vec{B}})$,
respectively.

\vspace{.3cm}
\begin{algorithm}[H]
\SetAlgoLined
\KwIn{$ \mathcal{A} \in \mathbb{R}^{\ell \times m \times n},\; \mathcal{B} \in
\mathbb{R}^{\ell \times p \times n}$ such that 
$\mathcal{A}^T*\mathcal{M}^{-1}*\mathcal{\vec{B}}\neq\mathcal{\vec{O}}$, 
$\mathcal{L} \in \mathbb{R}^{m \times m \times n}$ and 
$\mathcal{M} \in \mathbb{R}^{\ell \times \ell \times n}$ are SPD, $k\geq 1$}
$\beta_1 \leftarrow \|\mathcal{\vec{B}}\|_{\mathcal{M}^{-1}}$, $\mathcal{\vec{Q}}_1 \leftarrow
\frac{1}{\beta_1}\mathcal{\vec{B}}$\\
$\alpha_1 \leftarrow \|\mathcal{A}^T*\mathcal{M}^{-1}*\mathcal{\vec{Q}}_1\|_\mathcal{L}$, 
$\mathcal{\vec{W}}_1 \leftarrow \frac{1}{\alpha_1}(\mathcal{A}^T*\mathcal{M}^{-1}*\mathcal{\vec{Q}}_1)$\\
\For {$j = 1,2, \dots, k$}{
$\mathcal{\vec{Q}} \leftarrow \mathcal{A*L}*\mathcal{\vec{W}}_j - \alpha_j\mathcal{\vec{Q}}_j$\\
$\beta_{j+1} \leftarrow \|\mathcal{\vec{Q}}\|_{\mathcal{M}^{-1}}$, 
$\mathcal{\vec{Q}}_{j+1} \leftarrow \mathcal{\vec{Q}}/\beta_{j+1}$\\
\If {$j<k$}{
$\mathcal{\vec{W}} \leftarrow \mathcal{A}^T*\mathcal{M}^{-1}*\mathcal{\vec{Q}}_{j+1} - 
\beta_{j+1}\mathcal{\vec{W}}_j$\\
$\alpha_{j+1} \leftarrow \|\mathcal{\vec{W}}\|_\mathcal{L}$,
$\mathcal{\vec{W}}_{j+1} \leftarrow \mathcal{\vec{W}}/\alpha_{j+1}$\\
}
}
\caption{Partial weighted global tGKB (WG-tGKB) process}
\label{Alg: 11}
\end{algorithm}\vspace{.3cm}

Let $\mathcal{\vec{Y}}=\mathcal{W}_k \circledast y$. Then following a similar approach as
in Subsection \ref{sec4.1}, we reduce \eqref{4w} to the Tikhonov minimization problem in 
standard form
\begin{equation}\label{J1}
\min_{z\in \mathbb{R}^k}\{ \|\bar{B}_k z -
e_1 \beta_1\|^2_2 + \mu^{-1}\| z\|^2_2\}.
\end{equation}
The solution of \eqref{J1} determines an approximate solution of \eqref{4w}. We refer to 
this approach of computing an approximate solution of \eqref{4w} as the WG-tGKT method. 

The solution method for \eqref{4m} by applying the WG-tGKT method $p$ times to the
problems \eqref{33n} independently is referred to as the WG-tGKT$_p$ method. This 
$p$-method is described by Algorithm \ref{Alg: 13}.

\vspace{.3cm}
\begin{algorithm}[H]
\SetAlgoLined
\KwIn{$\mathcal{A}$, $\mathcal{\vec{B}}_1, \mathcal{\vec{B}}_2,\dots,\mathcal{\vec{B}}_p$,
$\mathcal{L}$, $\mathcal{M}$, $\delta_1, \delta_2, \dots, \delta_p$, $\eta > 1$, $ k_\text{init}= 2$}
\For{$j = 1,2,\dots, p$}{
$k \leftarrow k_\text{init}$, $\beta_1 \leftarrow \|\mathcal{\vec{B}}_j\|_{\mathcal{M}^{-1}}$,
$\mathcal{\vec{Q}}_1 \leftarrow \frac{1}{\beta_1}\mathcal{\vec{B}}_j$\\
Compute $\mathcal{W}_k$, $\mathcal{Q}_{k+1}$, and $\bar{B}_k$ by Algorithm \ref{Alg: 11}\\
Solve the minimization problem
\[
\min_{z \in \mathbb{R}^k} \| \bar{B}_kz - e_1 \beta_1 \|_2
\]
for $z_k$

\While{$\|\bar{B}_k z_k - e_1 \beta_1\|_2 \geq \eta \delta_j$}{
$k \leftarrow k+1$\\
$\mathtt{Go \; to \; step}$ 3
}
Determine the regularization parameter by the discrepancy principle, i.e., compute the zero $\mu_k>0$ of
\[
\varphi_k(\mu_k):=\|\bar{B}_k z_{j,\mu_k} - e_1 \beta_1\|_2^2-\eta^2 \delta_j^2
\]
and the associated solution $z_{j,\mu_k}$ of
\[
\min_{z \in \mathbb{R}^k}\left\| \left[ \begin{array}{c}
\bar{B}_k \\
\mu_k^{-1/2} I
\end{array}\right] z - \left[\begin{array}{c}
e_1 \beta_1\\
0
\end{array}\right] \right\|_2
\] \\
Compute $ \mathcal{Y}_{j,\mu_k} \leftarrow \mathcal{W}_k \circledast z_{j,\mu_k}$,
$\mathcal{\vec{X}}_{j,\mu_k} \leftarrow \mathcal{L} * \mathcal{Y}_{j,\mu_k}$ \\
}
\caption{The WG-tGKT$_p$ method for the approximate solution of \eqref{4m}.}
\label{Alg: 13}
\end{algorithm}\vspace{.3cm}

\section{Numerical Examples}\label{sec5}
All computations reported in this section are carried out in MATLAB R2021a with about $15$
significant decimal digits on a Dell computer running Windows 10 with 11th Gen Intel(R) 
Core(TM) i7-1165G7 $@$ 2.80GHz and 16 GB RAM. We compare the performance of the W-tGKT and
WG-tGKT methods for the solution of \eqref{4w} and \eqref{4m}, and the WGG-tGKT method for 
solving the latter problem. The implementation of the described methods is carried out in 
the spatial domain in Examples \ref{E2w}-\ref{E3w} using \eqref{spatialw}. All examples 
are concerned with image or video restoration.

Let $\mathcal{\vec{X}}_\text{method}\in\mathbb{R}^{m \times 1\times n}$ denote the 
approximate solution of \eqref{4w} computed by a particular solution method, and let 
$\mathcal{\vec{X}}_\text{true}\in\mathbb{R}^{m \times 1\times n}$ stand for the desired 
solution (e.g., a blur- and noise-free image). To compare the performance of the solution 
methods discussed, we tabulate the relative errors in the Frobenius norm,
\begin{equation}\label{relerr}
E_\text{method}=\frac{\|\mathcal{\vec{X}}_\text{method}-\mathcal{\vec{X}}_\text{true}\|_F}
{\|\mathcal{\vec{X}}_\text{true}\|_F}.
\end{equation}
We also display the Peak Signal-to-Noise Ratio,
\[
{\rm PSNR}:=10\,\log_{10}\bigg( \frac{{\rm MAX}_{\mathcal{\vec{X}}_{\rm true}}}
{\sqrt{\rm MSE}}\bigg),
\]
where ${\rm MAX}_{\mathcal{\vec{X}}_{\rm true}}$ stands for the maximum pixel value of the
desired image $\mathcal{\vec{X}}_{\rm true}$. The mean square error is given by
\[
{\rm MSE} = \frac{1}{mn} \sum_{i=1}^m \sum_{k=1}^n
\big(\mathcal{\vec{X}}_{\rm true} (i,1,k) - \mathcal{X}_{\rm method}(i,1,k)\big)^2.
\]

The relative error and the PSNR of computed approximate solutions 
$\mathcal{X}_\text{method}$ of \eqref{4m}, whose data is given by a tensor 
$\mathcal{B}\in\mathbb{R}^{m\times p\times n}$ for some $p>1$, are determined 
analogously. For the sake of brevity, we only explicitly discuss the situation when 
the data is a tensor slice $\mathcal{\vec{B}}\in\mathbb{R}^{m\times 1\times n}$.  

We generate a ``noise'' tensor $\mathcal{\vec{E}}$ that simulates the error in the data 
tensor $\mathcal{\vec{B}}=\mathcal{\vec{B}}_{\text{true}}+\mathcal{\vec{E}}$ in \eqref{4w}
and has a specified covariance tensor as described below.

\begin{prop}\label{prop5.1}
Let the entries of the tensor $\mathcal{\vec{E}}_\rho\in\mathbb{R}^{\ell\times 1\times n}$ be 
normally distributed with mean zero and variance $\rho>0$, and let $\mathcal{C}$ be a 
tensor of compatible size. Then the covariance tensor for 
$\mathcal{\vec{E}}=\mathcal{C}*\mathcal{\vec{E}}_\rho$ is $\rho^2\mathcal{C}*\mathcal{C}^T$. 
\end{prop}
\noindent
{\it Proof:} This result is well known in the situation when $\mathcal{\vec{E}}_\rho$ is a 
vector and $\mathcal{C}$ is a matrix; see, e.g., \cite[Lemma 2.1.1]{Bj} for a proof for
this case. The statistical properties of ${\mathcal{\vec{E}}}$ are independent of the
operations $\mathtt{fold}$ and $\mathtt{unfold}$ in the definition \eqref{tenA} of the 
t-product. We apply \cite[Lemma 2.1.1]{Bj} since $\mathtt{bcirc}(\mathcal{C})$ is a matrix
and ${\tt unfold}(\mathcal{\vec{E}}_\rho)$ is a vector. We also use the fact that 
$\mathtt{bcirc}(\mathcal{C}^T)=(\mathtt{bcirc}(\mathcal{C}))^T$, and 
$\mathtt{bcirc}(\mathcal{C}*\mathcal{C}^T) =
\mathtt{bcirc}(\mathcal{C})\cdot\mathtt{bcirc}(\mathcal{C}^T)$. Let $\mathtt{Cov}(\cdot)$ 
denote the covariance of a known quantity. Then 
\begin{eqnarray*}
\mathtt{Cov}(\mathcal{\vec{E}}) &=& {\tt fold}({\tt Cov}({\tt bcirc}(\mathcal{C})~
{\tt unfold}(\mathcal{\vec{E}}_\rho))) \\ 
&=& \rho^2{\tt fold}({\tt bcirc}(\mathcal{C}) \cdot({\tt bcirc}(\mathcal{C}))^T)\\ 
&=& \rho^2{\tt fold}({\tt bcirc}(\mathcal{C}*\mathcal{C}^T)\cdot{\tt unfold}(\mathcal{I}))
\\ &=& \rho^2\mathcal{C}*\mathcal{C}^T.  \hspace{5cm}\Box
\end{eqnarray*}

Let the tensor $\mathcal{\vec{E}}_\rho\in\mathbb{R}^{\ell\times 1 \times n}$ be defined as
in Proposition \ref{prop5.1} and let the tensor $\mathcal{C}$ be of compatible size. The 
entries of the tensor
\begin{equation}\label{Eten}
\mathcal{\vec{E}}:=\mathcal{C}*\mathcal{\vec{E}}_\rho
\end{equation}
simulate the noise with covariance tensor $\rho^2\mathcal{C}*\mathcal{C}^T$. We will use 
$\mathcal{\vec{E}}$ to simulate noise in 
$\mathcal{\vec{B}}\in\mathbb{R}^{\ell\times 1\times n}$, i.e., 
$\mathcal{\vec{B}}=\mathcal{\vec{B}}_{\text{true}}+\mathcal{\vec{E}}$, cf. \eqref{Btrue}. 
Introduce the scaled covariance tensor for $\mathcal{\vec{E}}$,
\begin{equation}\label{scaledcovten}
\mathcal{M}=\mathcal{C}*\mathcal{C}^T,
\end{equation}
where we assume that the tensor \eqref{scaledcovten} is positive definite. This is the 
generic situation. We refer to the quotient 
\begin{equation}\label{delta}
\widetilde{\delta}:=\frac{\|\mathcal{\vec{E}}\|_{\mathcal{M}^{-1}}}
{\|\mathcal{\vec{B}}_{\text{true}}\|_F} =\frac{\|\mathcal{\vec{E}}_\rho\|_F}
{\|\mathcal{\vec{B}}_{\text{true}}\|_F}
\end{equation}
as the \emph{noise level} of the error $\mathcal{\vec{E}}$ in $\mathcal{\vec{B}}$. 
Equality of the right-hand side of \eqref{delta} holds since 
\begin{eqnarray*}
\|\mathcal{\vec{E}}\|_{\mathcal{M}^{-1}}^2&=&
\|\mathcal{C}*\mathcal{\vec{E}}_\rho\|^2_{\mathcal{M}^{-1}}= 
\left(\mathcal{\vec{E}}_\rho^T*\mathcal{C}^T*\mathcal{M}^{-1}*
\mathcal{C}*\mathcal{\vec{E}}_\rho\right)_{(:,:,1)} \\
 &=& \left(\mathcal{\vec{E}}_\rho^T*
\mathcal{\vec{E}}_\rho\right)_{(:,:,1)} = \|\mathcal{\vec{E}}_\rho\|_F^2. 
\end{eqnarray*}

In the computed examples below, we prescribe the noise level \eqref{delta} and adjust 
$\rho$ in the noise tensor $\mathcal{\vec{E}}_\rho$ to obtain a tensor \eqref{Eten} that 
corresponds to the desired noise level. Specifically, since 
$\|\mathcal{\vec{E}}_\rho\|_F=\rho\|\mathcal{\vec{E}}_1\|_F$, we can adjust $\rho>0$ to 
obtain a ``noise tensor'' $\mathcal{\vec{E}}_\rho$ of a specified noise level. Knowledge 
of the noise level allows us to apply the discrepancy principle to determine the 
regularization parameter(s) as described in the previous sections.

In actual applications, the tensor $\mathcal{M}$ or an estimate thereof often are known 
and typically are symmetric positive definite. In our numerical experiments, we let
\begin{equation}\label{L1tensor}
\mathcal{M}=\mathcal{\widetilde{L}}_1^T*\mathcal{\widetilde{L}}_1+\omega\mathcal{I},
\end{equation}
with $\omega>0$. Here $\mathcal{\widetilde{L}}_1$ is a tensor, whose first frontal slice 
is the upper bidiagonal matrix
\[
\mathcal{\widetilde{L}}^{(1)}_1 =\frac{1}{2}\begin{bmatrix}
1 & -1 \\
& \ddots & \ddots \\
& & 1 & -1\\
& & & 1
\end{bmatrix}\in\mathbb{R}^{\ell \times \ell},
\]
and the remaining frontal slices 
$\mathcal{\widetilde{L}}_1^{(i)}\in\mathbb{R}^{\ell\times\ell}$, $i=2,3,\ldots,n$, are 
zero matrices. We compute the tensor Cholesky factorization of \eqref{L1tensor} by 
Algorithm \ref{cholesky} to obtain the Cholesky factor $\mathcal{R}$ and generate the
noise tensor by
\[
\mathcal{\vec{E}}=\mathcal{R}^T*\mathcal{\vec{E}}_\rho.
\]
Thus, equation \eqref{scaledcovten} holds with $\mathcal{C}=\mathcal{R}^T$. We adjust
$\rho>0$ to achieve a specified noise level as described above. 

The SPD tensors $\mathcal{D}_\gamma$ defined next are used as regularization tensor 
$\mathcal{L}$ in the Tikhonov minimization problems \eqref{4w} and \eqref{4m}. Let
\begin{equation}\label{tensorL}
\mathcal{D}_\gamma=\frac{1}{4}(\mathcal{\widetilde{L}}_2^T*\mathcal{\widetilde{L}}_2 + 
\alpha \mathcal{I}),  ~~~~\gamma \in \{1, 2\},
\end{equation}
where $\alpha>0$ and the tensor $\mathcal{\widetilde{L}}_2\in\mathbb{R}^{m\times m\times n}$ has the 
tridiagonal matrix
\[
{\mathcal{\widetilde{L}}}_2^{(1)} = \begin{bmatrix}
\gamma & -1 \\
-1 & 2 & -1 \\
& \ddots & \ddots & \ddots \\
& & -1 & 2 & -1\\
& & &-1 &\gamma
\end{bmatrix}\in\mathbb{R}^{m \times m},
\]
as its first frontal slice, and the remaining frontal slices 
${\mathcal{\widetilde{L}}_2}^{(i)}\in \mathbb{R}^{m \times m}$, $i=2,3,\ldots,n$, are zero
matrices. We do not require Cholesky factorization of the tensor \eqref{tensorL}.

The quality of restorations obtained for $\mathcal{L}=\mathcal{I}$ and 
$\mathcal{L}=\mathcal{D}_\gamma$ in the computed examples is compared. We use the 
parameters $\omega = 0.2$ in \eqref{L1tensor} and $\alpha=3$ in \eqref{tensorL} in most 
examples.  The influence on the computations of these parameters is discussed in Example 
\ref{E2w}. In all but the last example, we use $\gamma=1$ in \eqref{tensorL}. The last 
example compares this choice to $\gamma=2$. 

The symbol ``-'' in the tables indicates that the solution method solves several 
subproblems and carries out different numbers of bidiagonalization steps for the 
subproblems or uses several regularization parameters that take on different values. 
We determine the regularization parameter(s) by the bisection method over a specified 
interval using the discrepancy principle with safety parameter $\eta$. Specifically, we 
use the interval $[10^1,10^8]$ and $\eta = 1.01$ in Example \ref{E1w}, and the interval
$[10^1,10^7]$ and $\eta = 1.1$ in Examples \ref{E2w} and \ref{E3w}.

In all examples, the tensor $\mathcal{A}$ is a blurring operator that is constructed by 
using the function ${\tt blur}$ from \cite{Haa}. We use the {\tt squeeze} and {\tt twist} 
operators described in \cite{KBHH}, and the {\tt multi$\_$squeeze} and 
{\tt multi$\_$twist} operators defined in \cite{RU1} to store tensors in the desired 
format.

\begin{Ex}\label{E2w}(Color image restoration.)
This example illustrates the performance of the W-tGKT$_p$, WG-tGKT$_p$, and WGG-tGKT 
methods when applied to the restoration of the {\tt peppers} image of size 
$300\times 300\times 3$ shown in Figure \ref{Fig: 3w} (left). The image is stored in the 
RGB format, with each frontal slice corresponding to a color. The blurring tensor 
$\mathcal{A}\in \mathbb{R}^{300\times 300 \times 300}$ is generated by the MATLAB commands
\begin{equation}\label{zw}
z=\mathtt{[exp(-([0:band-1].^2)/(2\sigma^2)),zeros(1,N-band)]},
\end{equation}
\[
A = \frac{1}{\sigma\sqrt{2\pi}} \mathtt{toeplitz}(z), \;\; \mathcal{A}^{(i)} = A(i,1)A,
\;\; i = 1,2, \dots, 300, ~~ \sigma=3,~~ {\rm and} ~~\mathtt{band}=12.
\]

Let the true {\tt peppers} image (assumed to be unknown) be stored as 
$\mathcal{X}_{\rm true}\in \mathbb{R}^{300\times 3\times 300}$ by using the 
{\tt multi$\_$twist} operator. The blurred and noisy {\tt peppers} image is generated as
$\mathcal{B} = \mathcal{A}*\mathcal{X}_{\rm true} + \mathcal{E}$, where $\mathcal{E}$ is 
the noise tensor \eqref{Eten}. This image is displayed in Figure \ref{Fig: 3w} (right) for
$\widetilde{\delta}=10^{-3}$ by using the {\tt multi$\_$squeeze} operator.

\begin{figure}[!htb]
\hspace{1cm}
\minipage{0.42\textwidth}
\includegraphics[width=\linewidth]{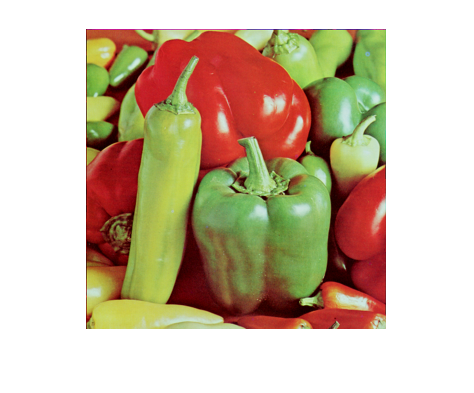} 
\endminipage\hfill \hspace{-2cm}
\minipage{0.42\textwidth}
\includegraphics[width=\linewidth]{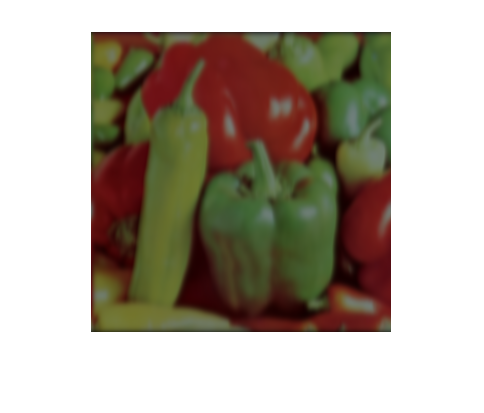}
\endminipage\hfill \vspace{-.9cm}
\caption{\small True image (left), and blurred and noisy image (right) for 
$\widetilde{\delta}=10^{-3}$.}\label{Fig: 3w}
\end{figure}

The dependence of the performance of the tGKT$_p$, WG-tGKT$_p$, and WGG-tGKT methods on 
the choice of $\omega>0$ for $\alpha=3$, and on the choice of $\alpha>0$ for $\omega=0.2$,
is illustrated by Figures \ref{Fig: exw1}-\ref{Fig: exw2} for 
$\widetilde{\delta}=10^{-2}$. We see from Figure \ref{Fig: exw1} that smaller values of 
$\omega$ often give higher PSNR values of the restored image and result in higher CPU
time requirement, whereas Figure \ref{Fig: exw2} shows that smaller $\alpha$ values often 
result in smaller PSNR and faster computations. These observations inform our choices of 
$\omega$ and $\alpha$ in all computed examples. Recall that the parameter $\omega$ affects 
the properties of the noise tensor $\mathcal{E}$, while $\alpha$ determines the 
regularization tensor $\mathcal{L}=\mathcal{D}_1$.

\begin{figure}[!htb]
\hspace{1cm}
\minipage{0.45\textwidth}
\includegraphics[width=\linewidth]{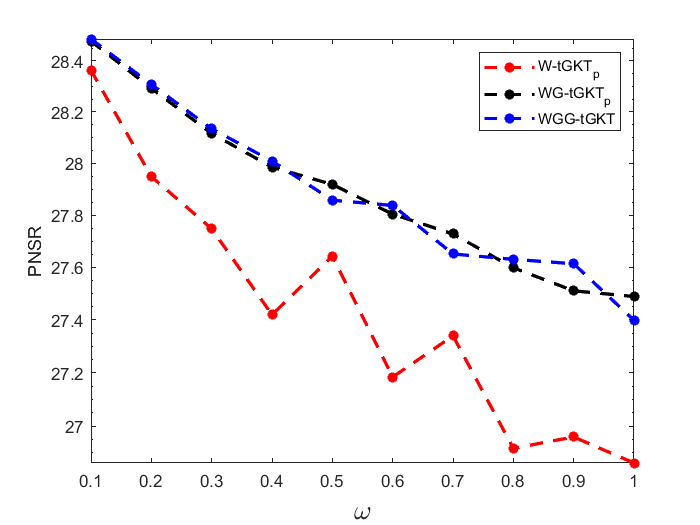} 
\endminipage\hfill \hspace{-1.5cm}
\minipage{0.45\textwidth}
\includegraphics[width=\linewidth]{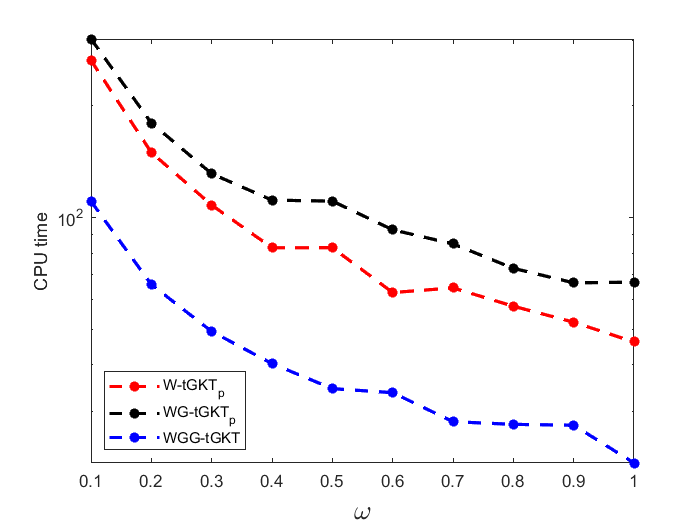}
\endminipage\hfill \vspace{-.2cm}
\caption{PSNR (left) and CPU time (right) for different $\omega$ values.}
\label{Fig: exw1}
\end{figure}

\begin{figure}[!htb]
\hspace{1cm}
\minipage{0.45\textwidth}
\includegraphics[width=\linewidth]{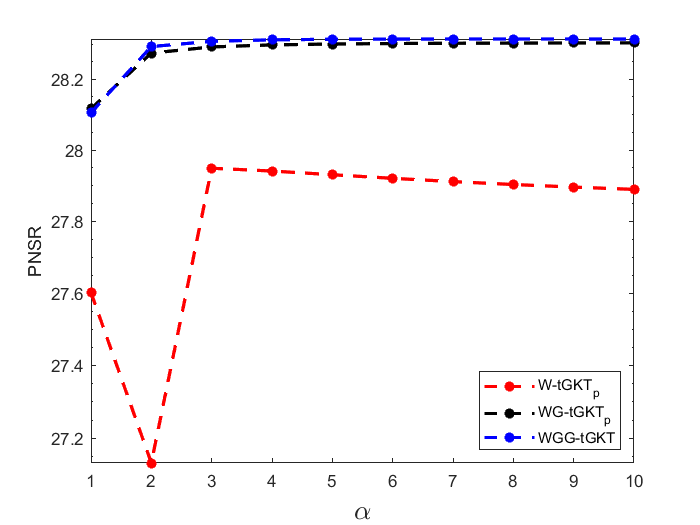} 
\endminipage\hfill \hspace{-1.5cm}
\minipage{0.45\textwidth}
\includegraphics[width=\linewidth]{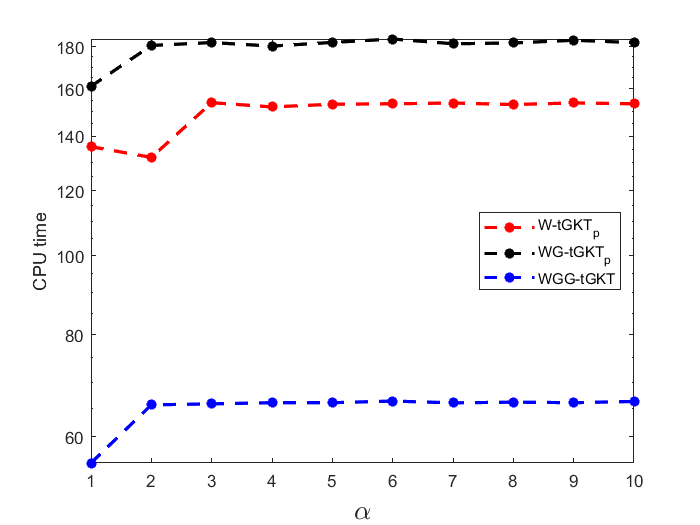}
\endminipage\hfill \vspace{-.2cm}
\caption{\small PSNR (left) and CPU time (right) for different $\alpha$ values.}
\label{Fig: exw2}
\end{figure}

\begin{figure}[!htb]
\hspace{1cm}
\minipage{0.42\textwidth}
\includegraphics[width=\linewidth]{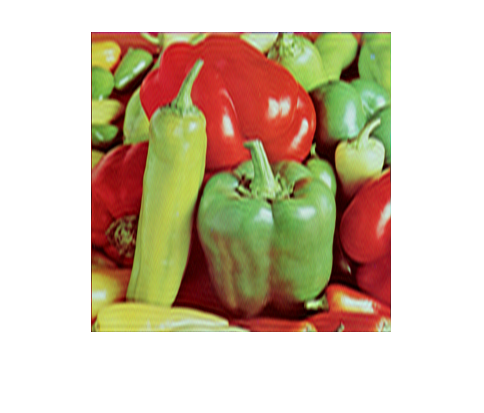} 
\endminipage\hfill \hspace{-2cm}
\minipage{0.42\textwidth}
\includegraphics[width=\linewidth]{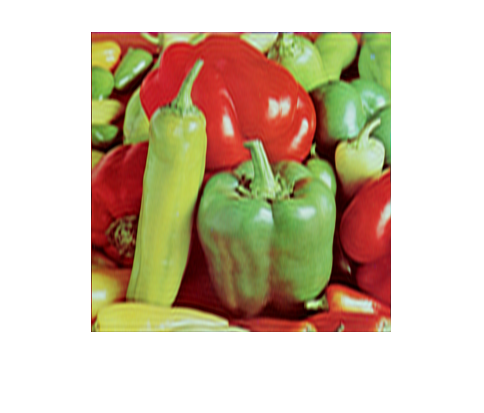}
\endminipage\hfill \vspace{-.9cm}
\caption{\small Reconstructed images by the W-tGKT$_p$ method (left) and the WGG-tGKT 
method after $67$ iterations (right) for $\widetilde{\delta}=10^{-3}$.}
\label{Fig: 4w}
\end{figure}

Figure \ref{Fig: 4w} shows restored images determined by the W-tGKT$_p$ and WGG-tGKT 
methods for $\mathcal{L}$ defined by \eqref{tensorL} with $\gamma=1$ and 
$\widetilde{\delta}=10^{-3}$, and Table \ref{Tab: 5.2} shows PSNR values and relative 
errors for each method as well as the CPU time required for the computations. Table 
\ref{Tab: 5.2} compares the performance of the W-tGKT$_p$, WG-tGKT$_p$, and WGG-tGKT
methods when applied to \eqref{4m} and to the minimization problem 
\be \label{4mspe}
\min_{\mathcal{X}\in\mathbb{R}^{m\times p\times n}}
\{\|\mathcal{A*X}-\mathcal{B}\|^2_F + 
\mu^{-1}\|\mathcal{X}\|^2_{\mathcal{L}^{-1}}\}.
\ee
We see from Table \ref{Tab: 5.2} that the use of weighted Frobenius norm in the 
fidelity term of \eqref{4m} is more appropriate when the error tensor $\mathcal{E}$ that 
simulates the noise in $\mathcal{B}$ is of the form \eqref{Eten} than the standard 
(unweighted) Frobenius norm in \eqref{4mspe}. Independently of the choice of $\mathcal{L}$
and noise level, the W-tGKT$_p$ method determines restorations of 
the worst quality for \eqref{4m}, and of the highest quality for \eqref{4mspe}. This 
method does not involve flattening and works separately with each lateral slice of the 
data tensor $\mathcal{B}$. For both choices of regularization tensor $\mathcal{L}$ in 
\eqref{4m}, the WG-tGKT$_p$ and WGG-tGKT methods give restorations of the highest quality 
for $\widetilde{\delta}=10^{-3}$ and $\widetilde{\delta}=10^{-2}$, respectively. Among all
methods considered for \eqref{4mspe}, the WGG-tGKT method results in restorations of the 
worst quality for both noise levels and for all choices of $\mathcal{L}$. This method 
requires less CPU time than the other methods for both minimization problems \eqref{4m} 
and \eqref{4mspe} because it works with the whole data tensor $\mathcal{B}$ at a time.

\begin{table}[h!]
\begin{center}
\begin{tabular}{cccccccc}
\multicolumn{8}{c}{The performance of the W-tGKT$_p$, WG-tGKT$_p$, and WGG-tGKT methods when applied to \eqref{4m}.}\\
\cmidrule(lr){1-8}
\multicolumn{1}{c}{$\widetilde{\delta}$} &\multicolumn{1}{c}{$\mathcal{L}$}&\multicolumn{1}{c}{Method}&\multicolumn{1}{c}{$k$} & \multicolumn{1}{c}{$\mu_k$}& \multicolumn{1}{c}{PSNR}&\multicolumn{1}{c}{Relative error}& \multicolumn{1}{c}{CPU time (secs)}
\\ \cmidrule(lr){1-8}
\multirow{6}{2em}{$10^{-3}$}&\multirow{3}{1em}{$\mathcal{D}_1$} & W-tGKT$_p$ &-&- &30.60&$5.5508\cdot 10^{-2}$ &$3.20\cdot 10^{3}$\\
&&WG-tGKT$_p$ &- &-&30.66& $5.5111\cdot 10^{-2}$ &$4.35\cdot 10^{3}$ \\
&&WGG-tGKT &67&$2.45\cdot 10^{4}$ &30.65& $5.5176\cdot 10^{-2}$ & $1.55\cdot 10^{3}$ \\ \cmidrule(lr){2-8}
&\multirow{3}{1em}{$\mathcal{I}$}& W-tGKT$_p$ &-&- &30.61&$5.5433\cdot 10^{-2}$ &$3.16\cdot 10^{3}$ \\
&&WG-tGKT$_p$ &- &-&30.64& $5.5262\cdot 10^{-2}$ &$5.12\cdot 10^{3}$ \\
&&WGG-tGKT &72&$2.41\cdot 10^{4}$ &30.62& $5.5356\cdot 10^{-2}$ & $1.76\cdot 10^{3}$ \\ \cmidrule(lr){1-8}
\multirow{6}{2em}{$10^{-2}$}&\multirow{3}{1em}{$\mathcal{D}_1$} & W-tGKT$_p$ &-&- &27.95&$7.5297\cdot 10^{-2}$ &$1.56\cdot 10^{2}$\\
&&WG-tGKT$_p$ &- &-&28.29& $7.2381\cdot 10^{-2}$ &$1.88\cdot 10^{2}$ \\
&&WGG-tGKT &14&$8.92\cdot 10^{2}$ &28.31& $7.2248\cdot 10^{-2}$ & $6.89\cdot 10^{1}$ \\ \cmidrule(lr){2-8}
&\multirow{3}{1em}{$\mathcal{I}$}& W-tGKT$_p$ &-&- &27.79&$7.6715\cdot 10^{-2}$ &$1.34\cdot 10^{2}$ \\
&&WG-tGKT$_p$ &- &-&28.31& $7.2265\cdot 10^{-2}$ &$1.86\cdot 10^{2}$ \\
&&WGG-tGKT &14&$9.39\cdot 10^{2}$ &28.31& $7.2214\cdot 10^{-2}$ & $6.85\cdot 10^{1}$ \\ \cmidrule(lr){1-8}
\multicolumn{8}{c}{The performance of the W-tGKT$_p$, WG-tGKT$_p$, and WGG-tGKT methods when applied to \eqref{4mspe}.}\\ \cmidrule(lr){1-8}
\multirow{6}{2em}{$10^{-3}$}&\multirow{3}{1em}{$\mathcal{D}_1$} & W-tGKT$_p$ &-&- &30.30&$5.7421\cdot 10^{-2}$ &$8.44\cdot 10^{2}$\\
&&WG-tGKT$_p$ &- &-&30.15& $5.8460\cdot 10^{-2}$ &$9.93\cdot 10^{2}$ \\
&&WGG-tGKT &32&$3.91\cdot 10^{4}$ &30.13& $5.8565\cdot 10^{-2}$ & $3.48\cdot 10^{2}$ \\ \cmidrule(lr){2-8}
&\multirow{3}{1em}{$\mathcal{I}$}& W-tGKT$_p$ &-&- &30.30&$5.7460\cdot 10^{-2}$ &$8.20\cdot 10^{2}$ \\
&&WG-tGKT$_p$ &- &-&30.13& $5.8601\cdot 10^{-2}$ &$1.14\cdot 10^{3}$ \\
&&WGG-tGKT &34&$3.98\cdot 10^{4}$ &30.11& $5.8723\cdot 10^{-2}$ & $4.03\cdot 10^{2}$ \\ \cmidrule(lr){1-8}
\multirow{6}{2em}{$10^{-2}$}&\multirow{3}{1em}{$\mathcal{D}_1$} & W-tGKT$_p$ &-&- &27.43&$7.9983\cdot 10^{-2}$ &$4.25\cdot 10^{1}$\\
&&WG-tGKT$_p$ &- &-&27.24& $8.1692\cdot 10^{-2}$ &$4.99\cdot 10^{1}$ \\
&&WGG-tGKT &7&$1.43\cdot 10^{3}$ &27.20& $8.2053\cdot 10^{-2}$ & $1.69\cdot 10^{1}$ \\ \cmidrule(lr){2-8}
&\multirow{3}{1em}{$\mathcal{I}$}& W-tGKT$_p$ &-&- &27.41&$8.0094\cdot 10^{-2}$ &$3.75\cdot 10^{1}$ \\
&&WG-tGKT$_p$ &- &-&27.32& $8.0918\cdot 10^{-2}$ &$5.39\cdot 10^{1}$ \\
&&WGG-tGKT &7&$1.22\cdot 10^{3}$ &27.21& $8.2004\cdot 10^{-2}$ & $1.69\cdot 10^{1}$ \\ \cmidrule(lr){1-8}
\end{tabular}
\end{center} \vspace{-.5cm}
\caption{\small Results for W-tGKT$_p$, WG-tGKT$_p$, and WGG-tGKT methods when applied to the restoration of {\tt peppers} image.}
\label{Tab: 5.2}
\end{table}

\end{Ex}

\begin{Ex}\label{E1w} (Medical imaging restoration.) 
We consider the restoration of an MRI image from MATLAB, and compare the weighted
Golub-Kahan-Tikhonov (W-GKT), W-tGKT, and WG-tGKT methods applied to the solution of 
\eqref{srhs}. The W-GKT method first reduces \eqref{srhs} to an equivalent problem 
involving a matrix and a vector, then applies the {\tt gen}$\_$GKB process described in 
\cite{CS} to determine an approximate solution of \eqref{srhs}. Our implementation of the 
W-GKT method is analogous to the implementations of the W-tGKT and WG-tGKT methods. The 
(unknown) true MRI image is shown on the left-hand side of Figure \ref{Fig: 1w}.

The frontal slices $\mathcal{A}^{(i)} \in \mathbb{R}^{256\times 256}$ of 
$\mathcal{A} \in \mathbb{R}^{256\times 256 \times 256}$ are generated by using a modified 
form of the function $\mathtt{blur}$ from \cite{Haa} with $\sigma = 4$ and ${\tt band}=7$.
Specifically, the blurring tensor $\mathcal{A}$ is generated with the MATLAB commands
\[
A = \frac{1}{\sigma\sqrt{2\pi}} \mathtt{toeplitz}([z(1)\;\mathtt{fliplr}(z(2:\mathtt{end}))],z), \;\; 
\mathcal{A}^{(i)} = A(i,1)A, \;\; i = 1,2, \dots, 256,
\]
where $z$ is defined in \eqref{zw}.

\begin{figure}[!htb]
\hspace{1cm}
\minipage{0.43\textwidth}
\includegraphics[width=\linewidth]{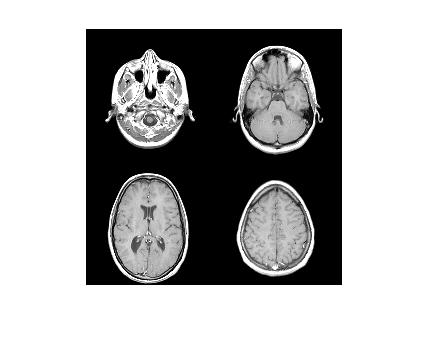}
\endminipage\hfill \hspace{-2cm}
\minipage{0.43\textwidth}
\includegraphics[width=\linewidth]{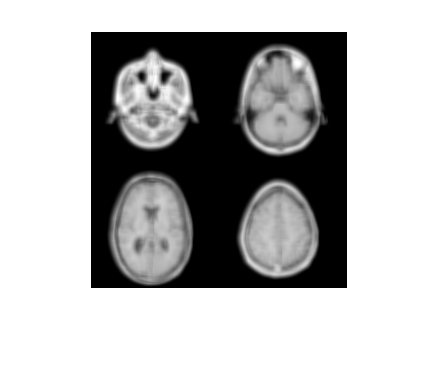}
\endminipage \hspace{1cm} \vspace{-.9cm}
\caption{\small True image (left), blurred noisy image (right) for $\widetilde{\delta}=10^{-3}$.}
\label{Fig: 1w}
\end{figure}

\begin{figure}[!htb]
\hspace{1cm}
\minipage{0.43\textwidth}
\includegraphics[width=\linewidth]{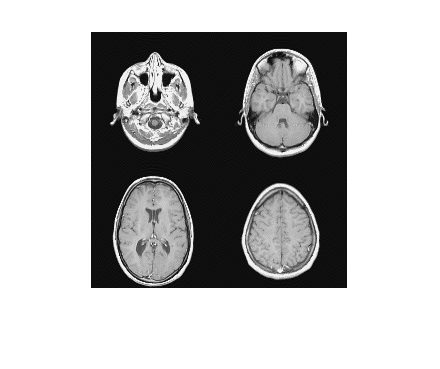} 
\endminipage\hfill \hspace{-2cm}
\minipage{0.43\textwidth}
\includegraphics[width=\linewidth]{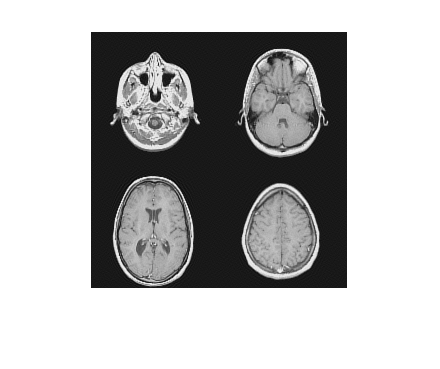}
\endminipage \hspace{1cm} \vspace{-.9cm}
\caption{\small Restored images determined by the W-tGKT method after $28$ iterations 
(left) and the WG-tGKT method after $107$ iterations (right) for 
$\widetilde{\delta}=10^{-3}$.}\label{Fig: 2w}
\end{figure}

We store the true MRI image of size $256 \times 256$ as a tensor 
$\mathcal{\vec{X}}_{\rm true} \in \mathbb{R}^{256\times 1 \times 256}$ by using the 
{\tt twist} operator. The blurred but noise-free image 
$\mathcal{\vec{B}}_{\rm true} \in \mathbb{R}^{256\times 1 \times 256}$ is generated by 
$\mathcal{\vec{B}}_{\rm true} = \mathcal{A}*\mathcal{\vec{X}}_{\rm true}$, and the blur- 
and noise-contaminated image 
$\mathcal{\vec{B}} = \mathcal{\vec{B}}_{\rm true} + \mathcal{\vec{E}}$ is displayed on the
right-hand side of Figure \ref{Fig: 1w} by using the {\tt squeeze} operator, where 
$\mathcal{\vec{E}}$ is a noise tensor analogous to \eqref{Eten}.

The restored images determined by the W-tGKT and WG-tGKT methods are displayed in Figure 
\ref{Fig: 2w} for the noise level $\widetilde{\delta}=10^{-3}$. Table \ref{Tab: 5.1} 
illustrates that the performance of the W-GKT, W-tGKT, and WG-tGKT methods depends on the 
choice of the tensor $\mathcal{L}$ and on the noise level. Independently of the choice of 
$\mathcal{L}$, the W-tGKT method, which does not involve flattening, yields restorations 
of the highest quality for $\widetilde{\delta}=10^{-3}$. The WG-tGKT method gives 
restorations of higher quality than the W-GKT method for $\widetilde{\delta}=10^{-3}$. The
latter method requires less CPU time than the former and is the fastest among the methods 
considered for both noise levels. The W-GKT and WG-tGKT methods yield restorations 
of almost the same quality for $\widetilde{\delta}=10^{-2}$, and require the same number 
of iterations for both noise levels independently of the choice of $\mathcal{L}$. These 
observations are based on the PSNR-values and relative errors shown in Table 
\ref{Tab: 5.1}. Regardless of the choice of $\mathcal{L}$ and the noise level, the W-tGKT 
method requires the least number of iterations, while the WG-tGKT method is the slowest. 

In summary, using the regularization tensor \eqref{tensorL} increases the quality of the
computed restorations slightly and reduces the number of iterations required by all methods
in our comparison. The method W-tGKT gives the most accurate restorations for the smaller
noise level. 

\begin{table}[h!]
\begin{center}
\begin{tabular}{cccccccc}
\cmidrule(lr){1-8}
\multicolumn{1}{c}{$\widetilde{\delta}$} &\multicolumn{1}{c}{$\mathcal{L}$}&\multicolumn{1}{c}{Method}&\multicolumn{1}{c}{$k$} & \multicolumn{1}{c}{$\mu_k$} &\multicolumn{1}{c}{PSNR}& \multicolumn{1}{c}{Relative error}& \multicolumn{1}{c}{CPU time (secs)}
\\ \cmidrule(lr){1-8}
\multirow{6}{2em}{$10^{-3}$}&\multirow{2}{1em}{$\mathcal{D}_1$} & W-GKT &107&$9.07\cdot 10^{5}$ &32.52&$5.4448\cdot 10^{-2}$ &$3.23\cdot 10^{2}$\\
&&W-tGKT &28&$1.85\cdot 10^{5}$&34.86& $4.1608\cdot 10^{-2}$ &$2.38\cdot 10^{2}$ \\
&&WG-tGKT &107&$2.44\cdot 10^{6}$&32.56& $5.4250\cdot 10^{-2}$ &$2.37\cdot 10^{3}$ \\ \cmidrule(lr){2-8}
&\multirow{3}{1em}{$\mathcal{I}$}& W-GKT &138&$3.85\cdot 10^{5}$ &32.19&$5.6571\cdot 10^{-2}$ &$3.20\cdot 10^{2}$ \\
&&W-tGKT &38&$2.60\cdot 10^{5}$&34.81& $4.1826\cdot 10^{-2}$ &$3.56\cdot 10^{2}$ \\
&&WG-tGKT &138&$3.97\cdot 10^{5}$&32.19& $5.6552\cdot 10^{-2}$ &$3.75\cdot 10^{3}$ \\ \cmidrule(lr){1-8}
\multirow{6}{2em}{$10^{-2}$}&\multirow{2}{1em}{$\mathcal{D}_1$} & W-GKT &25&$1.04\cdot 10^{4}$ &24.74&$1.3338\cdot 10^{-1}$ &$1.74\cdot 10^{1}$\\
&&W-tGKT &10&$3.10\cdot 10^{4}$&23.32& $1.5708\cdot 10^{-1}$ &$2.98\cdot 10^{1}$ \\
&&WG-tGKT &25&$1.07\cdot 10^{4}$&24.74& $1.3338\cdot 10^{-1}$ &$1.24\cdot 10^{2}$ \\ \cmidrule(lr){2-8}
&\multirow{3}{1em}{$\mathcal{I}$}& W-GKT &29&$8.74\cdot 10^{3}$ &24.74&$1.3347\cdot 10^{-1}$ &$1.46\cdot 10^{1}$ \\
&&W-tGKT &12&$3.87\cdot 10^{4}$&23.79& $1.4901\cdot 10^{-1}$ &$3.66\cdot 10^{1}$ \\
&&WG-tGKT &29&$8.85\cdot 10^{3}$&24.73& $1.3348\cdot 10^{-1}$ &$1.64\cdot 10^{2}$ \\ \cmidrule(lr){1-8}
\end{tabular}
\end{center} \vspace{-.5cm}
\caption{\small Results for the W-GKT, W-tGKT, and WG-tGKT methods when applied to the 
restoration of the {\tt MRI} image.}
\label{Tab: 5.1}
\end{table}
\end{Ex}

\begin{Ex}{(Video restoration.)}\label{E3w}
This example is concerned with the restoration of the first four consecutive frames of the
$\mathtt{Xylophone}$ video from MATLAB and compares the performance of the W-tGKT$_p$, 
WG-tGKT$_p$, and WGG-tGKT methods for the regularization tensors $\mathcal{L}=\mathcal{I}$ 
and for $\mathcal{L}$ defined by \eqref{tensorL} with $\gamma=1$ or $\gamma=2$. Each video
frame is in $\mathtt{MP4}$ format and has $240 \times 240$ pixels. The blurring operator 
$\mathcal{A}\in\mathbb{R}^{240 \times 240 \times 240}$ is generated similarly as in 
Example \ref{E2w} with $\sigma=2.5$ and ${\tt band}=12$.

\begin{figure}[!htb]
\hspace{1cm}
\minipage{0.42\textwidth}
\includegraphics[width=\linewidth]{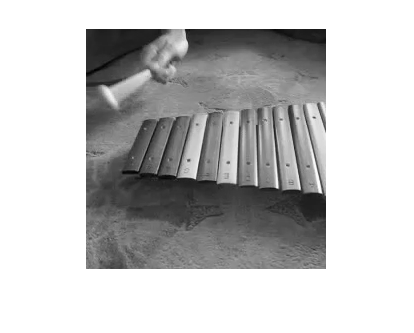} 
\endminipage\hfill \hspace{-2cm}
\minipage{0.42\textwidth}
\includegraphics[width=\linewidth]{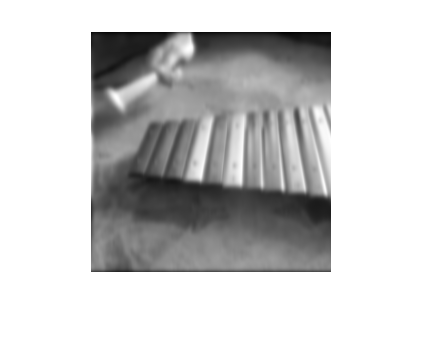}
\endminipage\hfill \vspace{-.9cm}
\caption{\small True fourth {\tt video} frame (left), and blurred and noisy fourth 
{\tt video} frame (right) for $\widetilde{\delta}=10^{-3}$.}
\label{Fig: 5w}
\end{figure}

\begin{figure}[!htb]
\hspace{1cm}
\minipage{0.42\textwidth}
\includegraphics[width=\linewidth]{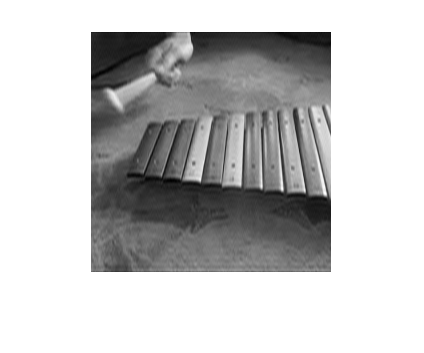} 
\endminipage\hfill \hspace{-2cm}
\minipage{0.42\textwidth}
\includegraphics[width=\linewidth]{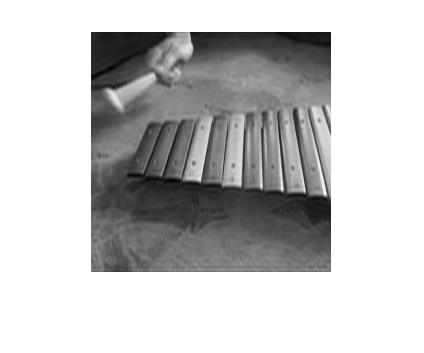}
\endminipage\hfill \vspace{-.9cm}
\caption{\small Restored fourth {\tt video} frame by the W-tGKT$_p$ method (left) and the 
WGG-tGKT method after $50$ iterations (right) for $\widetilde{\delta}=10^{-3}$.}
\label{Fig: 6w}
\end{figure}

The first four blur- and noise-free frames are stored as a tensor
$\mathcal{X}_\text{true}\in\mathbb{R}^{240 \times 4 \times 240}$ using the
$\mathtt{multi}\_\mathtt{twist}$ operator and are blurred by the tensor $\mathcal{A}$. We
generated the blur- and noise-contaminated frames as
$\mathcal{B} = \mathcal{A}*\mathcal{X}_{\rm true} + \mathcal{E}$, where $\mathcal{E}$ is
a noise tensor defined by \eqref{Eten}. The original fourth frame, and the 
corresponding blurred and noisy fourth frame are displayed in Figure \ref{Fig: 5w}
by using the {\tt squeeze} operator. The restored fourth frames determined by the 
W-tGKT$_p$ and WGG-tGKT methods are shown in Figure \ref{Fig: 6w} for $\widetilde{\delta}=10^{-3}$ and $\mathcal{L} = \mathcal{D}_1$. Comparing the PSNR values and relative errors, as well as the CPU times, 
displayed in Table \ref{Tab: 5.3}, the W-tGKT$_p$ method, independently of the choice of $\mathcal{L}$, yields restorations of the highest and of the least quality for $\widetilde{\delta}=10^{-3}$ and $\widetilde{\delta}=10^{-2}$, respectively. The WG-tGKT$_p$ method is seen to be the slowest and gives restorations of the highest quality for all choices of $\mathcal{L}$ for
$\widetilde{\delta}=10^{-2}$. The WGG-tGKT method is fastest for all choices of 
$\mathcal{L}$ and noise levels. Generally, the use of $\mathcal{L} = \mathcal{D}_2$ 
results in faster execution and gives restorations of higher quality for both noise levels 
than the other regularization tensors in our comparison.  

\begin{table}[h!]
\begin{center}
\begin{tabular}{cccccccc}
\cmidrule(lr){1-8}
\multicolumn{1}{c}{$\widetilde{\delta}$} &\multicolumn{1}{c}{$\mathcal{L}$}&\multicolumn{1}{c}{Method}&\multicolumn{1}{c}{$k$} & \multicolumn{1}{c}{$\mu_k$}& \multicolumn{1}{c}{PSNR}&\multicolumn{1}{c}{Relative error}& \multicolumn{1}{c}{CPU time (secs)}
\\ \cmidrule(lr){1-8}
\multirow{9}{2em}{$10^{-3}$}&\multirow{3}{1em}{$\mathcal{D}_1$} & W-tGKT$_p$ &-&- &34.03&$4.1173\cdot 10^{-2}$ &$8.22\cdot 10^{2}$\\
&&WG-tGKT$_p$ &- &-&34.03& $4.1185\cdot 10^{-2}$ &$1.25\cdot 10^{3}$ \\
&&WGG-tGKT &50&$3.32\cdot 10^{4}$ &34.03& $4.1193\cdot 10^{-2}$ & $3.47\cdot 10^{2}$ \\ \cmidrule(lr){2-8}
&\multirow{3}{1em}{$\mathcal{I}$}& W-tGKT$_p$ &-&- &34.11&$4.0797\cdot 10^{-2}$ &$8.24\cdot 10^{2}$ \\
&&WG-tGKT$_p$ &- &-&34.11& $4.0829\cdot 10^{-2}$ &$1.42\cdot 10^{3}$ \\
&&WGG-tGKT &54&$2.12\cdot 10^{4}$ &34.11& $4.0835\cdot 10^{-2}$ & $4.01\cdot 10^{2}$ \\ \cmidrule(lr){2-8}
&\multirow{3}{1em}{$\mathcal{D}_2$}& W-tGKT$_p$ &-&- &34.19&$4.0441\cdot 10^{-2}$ &$6.16\cdot 10^{2}$ \\
&&WG-tGKT$_p$ &- &-&34.14& $4.0648\cdot 10^{-2}$ &$1.12\cdot 10^{3}$ \\
&&WGG-tGKT &48&$1.42\cdot 10^{4}$ &34.15& $4.0632\cdot 10^{-2}$ & $3.20\cdot 10^{2}$ \\ \cmidrule(lr){1-8}
\multirow{9}{2em}{$10^{-2}$}&\multirow{3}{1em}{$\mathcal{D}_1$} & W-tGKT$_p$ &-&- &31.01&$5.8289\cdot 10^{-2}$ &$3.95\cdot 10^{1}$\\
&&WG-tGKT$_p$ &- &-&31.38& $5.5858\cdot 10^{-2}$ &$7.12\cdot 10^{1}$ \\
&&WGG-tGKT &12&$7.44\cdot 10^{2}$ &31.38& $5.5860\cdot 10^{-2}$ & $2.01\cdot 10^{1}$ \\ \cmidrule(lr){2-8}
&\multirow{3}{1em}{$\mathcal{I}$} & W-tGKT$_p$ &-&- &31.00&$5.8410\cdot 10^{-2}$ &$3.68\cdot 10^{1}$\\
&&WG-tGKT$_p$ &- &-&31.47& $5.5285\cdot 10^{-2}$ &$7.01\cdot 10^{1}$ \\
&&WGG-tGKT &12&$8.71\cdot 10^{2}$ &31.47& $5.5290\cdot 10^{-2}$ & $1.98\cdot 10^{1}$ \\ \cmidrule(lr){2-8}
&\multirow{3}{1em}{$\mathcal{D}_2$} & W-tGKT$_p$ &-&- &31.52&$5.5006\cdot 10^{-2}$ &$4.14\cdot 10^{1}$\\
&&WG-tGKT$_p$ &- &-&31.70& $5.3875\cdot 10^{-2}$ &$7.51\cdot 10^{1}$ \\
&&WGG-tGKT &11&$3.65\cdot 10^{3}$ &31.64& $5.4267\cdot 10^{-2}$ & $1.98\cdot 10^{1}$ \\ \cmidrule(lr){1-8}
\end{tabular}
\end{center} \vspace{-.5cm}
\caption{\small Results for the W-tGKT$_p$, WG-tGKT$_p$, and WGG-tGKT methods when applied
to the restoration of gray-scale {\tt video} frames.}
\label{Tab: 5.3}
\end{table}
\end{Ex}

\section{Conclusion}\label{sec6}
This paper extends the generalized Golub-Kahan bidiagonalization process described in 
\cite{CS} for matrices to third order tensors using a t-product. This results in the 
weighted t-product Golub-Kahan bidiagonalization (W-tGKB) process. Global versions of 
the latter process also are considered, namely, the weighted global t-product 
Golub-Kahan bidiagonalization (WG-tGKB) and the weighted generalized global t-product 
Golub-Kahan bidiagonalization (WGG-tGKB) processes. The W-tGKB process does not involve
flattening, but the global methods do. Only a few steps of the bidiagonalization processes
are required to solve the weighted Tikhonov regularization problems of our examples. This 
is typical for many image and video restoration problems. The use of a regularization
tensor $\mathcal{L \neq I}$ often results in higher quality restorations than when
$\mathcal{L = I}$.

The weighted t-product Golub-Kahan-Tikhonov (W-tGKT) and weighted global t-product 
Golub-Kahan-Tikhonov (WG-tGKT) regularization methods for the approximate solution of 
\eqref{4w} are considered. These methods are based on the W-tGKB and WG-tGKB processes, 
respectively. Independently of the choices of $\mathcal{L}$ considered, the W-tGKT method,
which does not involve flattening, yields the best or near-best quality restorations for 
$0.1\%$ noise level.

The W-tGKT and WG-tGKT methods also are applied $p$ times to determine an approximate 
solution of \eqref{4m}. This leads to the W-tGKT$_p$ and WG-tGKT$_p$ methods. The weighted
generalized global t-product Golub-Kahan-Tikhonov (WGG-tGKT) method for \eqref{4m} is also
discussed. This method differs from the W-tGKT$_p$ and WG-tGKT$_p$ methods in that it uses
the WGG-tGKB process and works with a large amount of data at a time.

The WGG-tGKT method is the fastest, while the WG-tGKT$_p$ method is the slowest,
independently of the choice of $\mathcal{L}$ and noise levels. Both methods involve 
flattening since they require additional product definition to the t-product. Generally, 
working with one lateral slice of the data tensor at a time is seen to give restorations 
of higher quality than working with all lateral slices simultaneously.


\begin{thebibliography}{27}

\label{sec:bib}

\bibitem{A}
M. Oriel, Generalized Golub-Kahan bidiagonalization and stopping criteria, SIAM J. 
Matrix Anal. Appl., 34 (2013), pp. 571--592.

\bibitem{AO}
M. Arioli and D. Orban, Iterative methods for symmetric quasi-definite linear systems – 
Part I: Theory. Cahier du GERAD G-2013-32. Montr\'eal, Canada: GERAD, Montr\'eal, QC, 2013.

\bibitem{ABH}
S. R. Arridge, M. M. Betcke, and L. Harhanen, Iterated preconditioned LSQR method for 
inverse problems on unstructured grids, Inverse Problems, 30 (2014), Art. 075009.

\bibitem{B}
S. J. Benbow, Solving generalized least-squares problems with LSQR, SIAM J. Matrix Anal.
Appl., 21 (1999), pp. 166--177.

\bibitem{BJNR}
F. P. A. Beik, K. Jbilou, M. Najafi-Kalyani, and L. Reichel, Golub-Kahan bidiagonalization
for ill-conditioned tensor equations with applications, Numer. Algorithms, 84 (2020), pp.
1535--1563.

\bibitem{BNR}
F. P. A. Beik, M. Najafi-Kalyani, and L. Reichel, Iterative Tikhonov regularization of
tensor equations based on the Arnoldi process and some of its generalizations, Appl.
Numer. Math., 151 (2020), pp. 425--447.

\bibitem{BIJS}
F. P. A. Beik, A. El Ichi, K. Jbilou, and R. Sadaka, Tensor extrapolation methods with 
applications, Numer. Algorithms, 87 (2021), pp. 1421--1444. 

\bibitem{Bj}
\AA. Bj\"orck, Numerical Methods in Matrix Computations, Springer, Cham, 2015.

\bibitem{CS}
J. Chung and A. Saibaba, Generalized hybrid iterative methods for large-scale Bayesian 
inverse problems, SIAM J. Sci. Comput., 39 (2007), pp. S24--S46.

\bibitem{CSBW}
J. Chung, A. Saibaba, M. Brown, and E. Westman, Efficient generalized Golub-Kahan based 
methods for dynamic inverse problems. Inverse Problems, 34 (2018), Art. 024005.

\bibitem{GIJB}
M. El Guide, A. El Ichi, K. Jbilou, and F. P. A. Beik, Tensor GMRES and Golub-Kahan
bidiagonalization methods via the Einstein product with applications to image and video
processing, https://arxiv.org/pdf/2005.07458.pdf

\bibitem{GIJS}
M. El Guide, A. El Ichi, K. Jbilou, and R. Sadaka, Tensor Krylov subspace methods via the
T-product for color image processing, June 2020, https://arxiv.org/pdf/2006.07133.pdf

\bibitem{IGJ}
A. El Ichi, M. El Guide, and K. Jbilou, Discrete cosine transform LSQR and GMRES methods 
for multidimensional ill-posed problems, March 2021. https://arxiv.org/pdf/2103.11847.pdf

\bibitem{EHN}
H. W. Engl, M. Hanke, and A. Neubauer, Regularization of Inverse Problems, Kluwer,
Dordrecht, 1996.

\bibitem{FRR}
C. Fenu, L. Reichel, and G. Rodriguez, GCV for Tikhonov regularization via global
Golub-Kahan decomposition, Numer. Linear Algebra Appl., 23 (2016), pp. 467--484.

\bibitem{Haa}
P. C. Hansen, Regularization Tools, version 4.0, for MATLAB 7.3. Numer. Algorithms, 46
(2007), pp. 189--194.

\bibitem{Ha}
P. C. Hansen, Rank-Deficient and Discrete Ill-Posed Problems, SIAM, Philadelphia, 1998.

\bibitem{KBHH}
M. E. Kilmer, K. Braman, N. Hao, and R. C. Hoover, Third-order tensors as operators on
matrices: A theoretical and computational framework with applications in imaging, SIAM
J. Matrix Anal. Appl., 34 (2013), pp. 148--172.

\bibitem{KM}
M. E. Kilmer and C. D. Martin, Factorization strategies for third-order tensors, Linear
Algebra Appl., 434 (2011), pp. 641--658.

\bibitem{Ki}
S. Kindermann, Convergence analysis of minimization-based noise level-free parameter
choice rules for linear ill-posed problems, Electron. Trans. Numer. Anal., 38 (2011),
pp. 233--257.

\bibitem{KR}
S. Kindermann and K. Raik, A simplified L-curve method as error estimator, Electron.
Trans. Numer. Anal., 53 (2020), pp. 217--238.

\bibitem{KB}
T. G. Kolda and B. W. Bader, Tensor decompositions and applications, SIAM Rev., 51 (2009),
pp. 455--500.

\bibitem{OA}
D. Orban and M. Arioli, Iterative Solution of Symmetric Quasi-Definite Linear Systems,
SIAM, Philadelphia, 2017.

\bibitem{SCP}
A. K. Saibaba, J. Chung, and K. Petroske, Efficient Krylov subspace methods for 
uncertainty quantification in large Bayesian linear inverse problems, Numer. Linear 
Algebra Appl., 27 (2020) Art. e2325. 

\bibitem{RR}
L. Reichel and G. Rodriguez, Old and new parameter choice rules for discrete ill-posed
problems, Numer. Algorithms, 63 (2013), pp. 65--87.

\bibitem{RU1}
L. Reichel and U. O. Ugwu, Tensor Golub-Kahan-Tikhonov methods applied to the solution of
ill-posed problem with a t-product structure, Numer. Linear Algebra Appl., in press.

\bibitem{RU2}
L. Reichel and U. O. Ugwu, Tensor Arnoldi-Tikhonov and GMRES-type methods for ill-posed
problem with t-product structure, 2020, submitted for publication.

\bibitem{RU3}
L. Reichel, and U. O. Ugwu, Tensor Krylov subspace methods with an invertible linear 
transform product applied to image processing, Appl. Numer. Math., 166 (2021), pp. 
186--207.


\end{thebibliography}
\end{document}